# Towards computable flows and robust estimates for inf-sup stable FEM applied to the time-dependent incompressible Navier–Stokes equations

**Philipp W. Schroeder** · **Christoph Lehrenfeld** ·
**Alexander Linke** · **Gert Lube**

**Abstract** Inf-sup stable FEM applied to time-dependent incompressible Navier–Stokes flows are considered. The focus lies on robust estimates for the kinetic and dissipation energies in a twofold sense. Firstly, pressure-robustness ensures the fulfilment of a fundamental invariance principle and velocity error estimates are not corrupted by the pressure approximability. Secondly, *Re*-semi-robustness means that constants appearing on the right-hand side of kinetic and dissipation energy error estimates (including Gronwall constants) do not explicitly depend on the Reynolds number. Such estimates rely on the essential regularity assumption $\nabla u \in L^1(0,T;L^\infty(\Omega))$ which is discussed in detail. In the sense of best practice, we review and establish pressure- and *Re*-semi-robust estimates for pointwise divergence-free $H^1$-conforming FEM (like Scott–Vogelius pairs or certain isogeometric based FEM) and pointwise divergence-free $H(\mathrm{div})$-conforming discontinuous Galerkin FEM. For convection-dominated problems, the latter naturally includes an upwind stabilisation for the velocity which is not gradient-based.

Alexander Linke
Weierstrass Institute
D-10117 Berlin, Germany
E-mail: alexander.linke@wias-berlin.de
ORCID: https://orcid.org/0000-0002-0165-2698

Philipp W. Schroeder & Christoph Lehrenfeld & Gert Lube
Institute for Numerical and Applied Mathematics
Georg-August-University Göttingen
D-37083 Göttingen, Germany
E-mail: {p.schroeder,lehrenfeld,lube}@math.uni-goettingen.de
ORCID P.W. Schroeder: https://orcid.org/0000-0001-7644-4693
ORCID C. Lehrenfeld: https://orcid.org/0000-0003-0170-8468

## 1 Introduction

We consider the time-dependent incompressible Navier–Stokes equations [66,61,27]

$$\begin{cases} \partial_t \boldsymbol{u} - \nu\Delta\boldsymbol{u} + (\boldsymbol{u}\cdot\nabla)\boldsymbol{u} + \nabla p = \boldsymbol{f} & \text{in } (0,T]\times\Omega, \quad (1a)\\ \nabla\cdot\boldsymbol{u} = 0 & \text{in } (0,T]\times\Omega, \quad (1b)\\ \boldsymbol{u} = \boldsymbol{0} & \text{on } [0,T]\times\partial\Omega, \quad (1c)\\ \boldsymbol{u}(0,\boldsymbol{x}) = \boldsymbol{u}_0(\boldsymbol{x}) & \text{for } \boldsymbol{x}\in\Omega. \quad (1d) \end{cases}$$

For the space dimension $d \in \{2,3\}$, $\Omega \subset \mathbb{R}^d$ denotes a connected bounded Lipschitz domain. Moreover, $\boldsymbol{u}\colon (0,T]\times\Omega \to \mathbb{R}^d$ indicates the velocity field, $p\colon (0,T]\times\Omega \to \mathbb{R}$ is the (zero-mean) kinematic pressure, $\boldsymbol{f}\colon (0,T]\times\Omega \to \mathbb{R}^d$ represents external body forces and $\boldsymbol{u}_0\colon \Omega \to \mathbb{R}^d$ stands for a suitable initial condition for the velocity. The underlying fluid is assumed to be Newtonian with constant (dimensionless) kinematic viscosity $0 < \nu \ll 1$.

There are references regarding the historical development of finite element methods (FEM) for the Navier–Stokes problem (1) until 2016; cf., for example, the monograph [39]. A summary of very recent results for $\boldsymbol{H}^1$-conforming FEM, together with several open problems, can be found in the review paper [40].

A relatively new aspect in the FE analysis applied to incompressible flows is ‘pressure-robustness’ [41]. In its most general form, pressure-robustness of a numerical method is defined by its ability to fulfil the following requirement: if the exact solution $\boldsymbol{u}$ of (1) belongs to the approximation space $\boldsymbol{V}_h$, i.e. if $\boldsymbol{u}\in\boldsymbol{V}_h$, then the discrete solution $\boldsymbol{u}_h$ coincides with the exact one, that is, $\boldsymbol{u}_h = \boldsymbol{u}$. In certain physical regimes of the incompressible Navier–Stokes equations — i.e., in certain benchmarks — pressure-robustness allows to use less expensive discretisation schemes without losing accuracy [49,1]. As a consequence, the following fundamental invariance principle transfers from the continuous level to the discretised case: Replacing the source term $\boldsymbol{f}$ by $\boldsymbol{f}+\nabla\psi$ changes the solution $(\boldsymbol{u},p)$ to $(\boldsymbol{u},p+\psi)$. For example, in a potential flow, $(\boldsymbol{u}\cdot\nabla)\boldsymbol{u}$ can be very large but it is a gradient and therefore balanced by the pressure gradient and thus does not have any impact on the velocity field. Only recently it has been shown that high Reynolds number potential flows are really challenging for the numerical solution with standard low-order Galerkin-FEM [50,41].

A well-known important consequence for methods which are not pressure-robust is that already for the steady incompressible Stokes problem the velocity error estimates for kinetic and dissipation energies are corrupted by the pressure approximability multiplied by $\nu^{-1/2}$ [41,49]. Note that the mechanism responsible for the excitation of this kind of numerical error is a completely linear phenomenon. Exactly divergence-free FEM are naturally pressure-robust, but classical inf-sup stable velocity-pressure pairs like Taylor–Hood FEM are usually not pressure-robust. In fact, such classical, inf-sup stable mixed finite elements that relax the divergence constraint are usually prone to the locking phenomenon of poor mass conservation [49,1]. Fortunately, recent research allows to slightly modify such methods in order to make them pressure-robust by so-called velocity reconstructions; for example for the Stokes problem, we refer to [47,48,41,43].

However, in this article, we focus on a different important aspect in the continuous-in-time numerical analysis; namely, the worst case behaviour of the velocity error due to the nonlinearity of the convection term in the time-dependent setting. This is reflected in the numerical error analysis by Gronwall constants depending at least exponentially on time. Indeed, in case of $0 < \nu \ll 1$, in many estimates available in the literature the constant $C$ in the exponential growth $\exp(Ct)$ in fact depends on the Reynolds number $Re$ (respectively, on $\nu^{-1}$) or even powers of $Re$; see Subsection 4.1. Obviously, such error estimates can describe a sensible error behaviour only for ultra short time intervals. The value of these estimates is that they predict correctly the convergence behaviour of the velocity errors with respect to space discretisation; although they involve huge constants in the estimates. In view of this situation, numerical analysts frequently argue that these error estimates might not be sharp. Following the original proposal by [60] for scalar diffusion-advection problems, error estimates where the constants appearing on the right-hand side (including Gronwall constants) do not explicitly depend on the Reynolds number are called '*Re*-semi-robust'.

Partially, the problems in the numerical analysis come from very weak assumptions on the exact solution $\boldsymbol{u}$ and the data. It turns out that error estimates can be improved considerably under the essential regularity condition $\nabla \boldsymbol{u} \in L^1(0,T;\boldsymbol{L}^\infty(\Omega))$. We summarise some physical implications of this stronger regularity condition in Subsection 2.2. On the other hand, the numerical analysis will show that for this class of flows the right-hand side of error estimates grows relatively mildly with $\exp(Ct)$, where $C$ is not explicitly dependent on the Reynolds number. Therefore, in this article, this class of flows is called 'computable'.

The stronger regularity condition has been used frequently in the literature, even in the limit case of incompressible Euler flow $\nu = 0$; cf., for example, the monograph [52] or the review [5]. In order to obtain *Re*-semi-robust error estimates for problem (1), [13] is presumably the first work which takes advantage of this regularity assumption in the analysis of a CIP-stabilised FEM with equal-order approximation of velocity and pressure. For an equal-order method with local projection stabilisation (LPS), we refer to the recent work [32]. However, using non-inf-sup stable methods excludes the possibility of obtaining pressure-robustness in [13,32]. Note that the concept of pressure-robustness goes beyond the question of optimal $h$-convergence rates, which are indeed proved in [13,32], since it is a robustness property of the velocity error. Concerning the debate on the optimality of classical mixed methods versus pressure-robust mixed methods, the reader is referred to [1, Section 4]. For $\boldsymbol{H}^1$-conforming inf-sup stable FEM, in [3] the combination with grad-div stabilisation in some different energy norm led to *Re*-semi-robust error estimates which were sharpened in [25]. The work in [33] deepens the results; in particular for optimal pressure estimates.

The main purpose of the present paper is to review the state-of-the-art concerning error estimates for exactly divergence-free FEM for problem (1). In particular, we concentrate on *Re*-semi-robustness and pressure-robustness for the velocity estimates. Due to the inherent pressure-robustness, it is possible to separate velocity and pressure completely in the error analysis. Therefore, we focus exclusively on velocity estimates. In such a setting, this paper offers a unified approach to continuous-in-time error estimates for exactly divergence-free $\boldsymbol{H}^1$-conforming and only $\boldsymbol{H}(\mathrm{div})$-conforming FEM. In particular, the extension to $\boldsymbol{H}(\mathrm{div})$-conforming FEM for problem (1) with $0 < \nu \ll 1$ is original. Results for such FEM in the case of the incompressible Euler equations with $\nu = 0$ can be found in [35,53].

In our opinion, using exactly divergence-free methods has the following main advantages: They are inherently pressure-robust. As shown in Section 5, one can obtain *Re*-semi-robust estimates without any additional stabilisation or skew-symmetrisation, thereby facil-itating the numerical analysis. Furthermore, the fact that less stabilisation is required allows to achieve discretisations where the amount of necessary numerical dissipation is minimised. In addition, the conservation properties of the exact solution to (1) (as for example conservation of mass, energy and momentum) are naturally transferred to the discrete solution. In particular, divergence-free methods have a healthy and clean energy balance *a priori*.

Concerning numerical experiments we show exemplarily that for a planar standing vortex problem (or periodic lattice flow), the exponential growth of the $Re$-semi-robust Gronwall-based error estimates can be observed at least qualitatively also in practice. Furthermore, for a time-dependent potential flow, we show that divergence-free methods are clearly superior to classical mixed FEM.

*Organisation of the article:* In Section 2 the continuous Navier–Stokes problem is briefly recalled and the meaningfulness of our essential regularity assumption is discussed. Afterwards, Section 3 lays the foundation for a unified treatment of FEM for the time-dependent Navier–Stokes problem. For classical $\boldsymbol{H}^1$-conforming methods, a brief treatise and recollection of Galerkin-FEM and grad-div stabilisation is provided in Section 4. Then, moving to exactly divergence-free FEM, Section 5 treats pressure- and $Re$-semi-robust error estimates for both $\boldsymbol{H}^1$- and $\boldsymbol{H}(\mathrm{div})$-conforming methods in a unified setting. Some numerical experiments are also conducted. After a brief survey about some open problems in Section 6, the main part of this work is concluded in Section 7. Computational aspects of $\boldsymbol{H}(\mathrm{div})$-conforming methods are addressed in the Appendix.

## 2 Continuous Navier–Stokes problem

*Notation:* In what follows, for $K \subseteq \Omega$ we use the standard Sobolev spaces $W^{m,p}(K)$ for scalar-valued functions with associated norms $\|\cdot\|_{W^{m,p}(K)}$ and seminorms $|\cdot|_{W^{m,p}(K)}$ for $m \geqslant 0$ and $p \geqslant 1$. Spaces and norms for vector- and tensor-valued functions are indicated with bold letters. We obtain the Lebesgue space $W^{0,p}(K) = L^p(K)$ and the Hilbert space $W^{m,2}(K) = H^m(K)$. Additionally, the closed subspaces $H_0^1(K)$ consisting of $H^1(K)$-functions with vanishing trace on $\partial K$ and the set $L_0^2(K)$ of $L^2(K)$-functions with zero mean in $K$ play an important role. The $L^2(K)$-inner product is denoted by $(\cdot,\cdot)_K$ and, if $K = \Omega$, we usually omit the domain completely when no confusion can arise. Furthermore, with regard to time-dependent problems, given a Banach space $\boldsymbol{X}$ and a time instance $t$, the Bochner space $L^p(0,t;\boldsymbol{X})$ for $p \in [1,\infty]$ is used. In the case $t = T$, we frequently use the abbreviation $L^p(\boldsymbol{X}) = L^p(0,T;\boldsymbol{X})$. The dual space of $\boldsymbol{X}$ is denoted by $\boldsymbol{X}^*$.

### 2.1 Continuous problem

With $\boldsymbol{V} = \boldsymbol{H}_0^1(\Omega)$ and $Q = L_0^2(\Omega)$, we introduce the space

$$\boldsymbol{V}^{\mathrm{div}} = \{\boldsymbol{v} \in \boldsymbol{V} : (q, \nabla\cdot\boldsymbol{v}) = 0,\ \forall q \in Q\} \tag{2}$$

of weakly divergence-free velocities. If $\boldsymbol{v} \in \boldsymbol{V}^{\mathrm{div}}$, then $\nabla\cdot\boldsymbol{v} = 0$ almost everywhere in $\Omega$. The largest space in which one can work comfortably with the divergence is

$$\boldsymbol{H}(\mathrm{div};\Omega) = \left\{\boldsymbol{v} \in \boldsymbol{L}^2(\Omega) :\ \nabla\cdot\boldsymbol{v} \in L^2(\Omega)\right\}. \tag{3}$$

We remark that the statement $\nabla\cdot\boldsymbol{v}\in L^2(\Omega)$ in this definition means that the *distributional divergence* of $\boldsymbol{v}$ lies in $L^2$. For the importance of the notion of the *distributional divergence* with respect to understanding the significance of pressure-robustness, we refer to [41]. Analogously to $\boldsymbol{V}^{\mathrm{div}}$, we define

$$\boldsymbol{H}^{\mathrm{div}} = \left\{\boldsymbol{v}\in\boldsymbol{H}(\mathrm{div};\Omega):\ \nabla\cdot\boldsymbol{v}=0,\ \boldsymbol{v}\cdot\boldsymbol{n}\big|_{\partial\Omega}=0\right\}, \tag{4}$$

where $\boldsymbol{n}$ denotes the outer unit normal vector to $\partial\Omega$. For the following error analysis, velocity and pressure solutions are assumed to belong to the spaces

$$\boldsymbol{V}^T = \left\{\boldsymbol{v}\in L^2(0,T;\boldsymbol{V}):\ \partial_t\boldsymbol{v}\in L^2\big(0,T;\boldsymbol{L}^2(\Omega)\big)\right\} \quad \text{and} \quad Q^T = L^2(0,T;Q). \tag{5}$$

Recently it has been shown for the solution $\boldsymbol{u}_s$ of the evolutionary Stokes problem with inhomogeneous Dirichlet data and $\boldsymbol{f}=\boldsymbol{0}$ that $\partial_t\boldsymbol{u}_s\in L^2\big(\boldsymbol{L}^2\big)$ indeed holds [16]. Thus, provided $\boldsymbol{f}\in L^2\big(\boldsymbol{L}^2\big)$, the following problem on the continuous level is obtained:

$$\begin{cases}\text{Find } (\boldsymbol{u},p)\in\boldsymbol{V}^T\times Q^T \text{ with } \boldsymbol{u}(0)=\boldsymbol{u}_0\in\boldsymbol{H}^{\mathrm{div}} \text{ s.t., } \forall\,(\boldsymbol{v},q)\in\boldsymbol{V}\times Q, & \text{(6a)}\\ (\partial_t\boldsymbol{u},\boldsymbol{v})+\nu a(\boldsymbol{u},\boldsymbol{v})+c(\boldsymbol{u};\boldsymbol{u},\boldsymbol{v})+b(\boldsymbol{v},p)-b(\boldsymbol{u},q)=(\boldsymbol{f},\boldsymbol{v}). & \text{(6b)}\end{cases}$$

Here, the multilinear forms are given by

$$a(\boldsymbol{w},\boldsymbol{v}) = \int_\Omega \nabla\boldsymbol{w}:\nabla\boldsymbol{v}\,\mathrm{d}\boldsymbol{x}, \qquad c(\boldsymbol{\beta};\boldsymbol{w},\boldsymbol{v}) = \int_\Omega (\boldsymbol{\beta}\cdot\nabla)\boldsymbol{w}\cdot\boldsymbol{v}\,\mathrm{d}\boldsymbol{x}, \tag{7a}$$

$$b(\boldsymbol{w},q) = -\int_\Omega q(\nabla\cdot\boldsymbol{w})\,\mathrm{d}\boldsymbol{x}. \tag{7b}$$

**Remark 2.1** Concerning the regularity of the forcing term, on the continuous level, the problem could be posed using the less restrictive assumption $\boldsymbol{f}\in L^2(\boldsymbol{V}^*)$. However, in Section 5 we also deal with discretisations which are not $\boldsymbol{H}^1$-conforming. In such a situation, rough right-hand sides lead to technical difficulties which we omit by assuming $\boldsymbol{f}\in L^2\big(\boldsymbol{L}^2\big)$; cf. [26, Remark 4.9]. Another problem with rough forcing terms, even for $\boldsymbol{H}^1$-conforming methods, is that energy estimates can generally not be expected to be independent of $\nu^{-1}$; cf. [64, Remark 3.2].

**Remark 2.2** The theory concerning existence and regularity of Navier–Stokes solutions gives the following result; cf. [6, 10, 39]. To (6b), there exists a weak solution

$$\boldsymbol{u}\in L^2\Big(0,T;\boldsymbol{V}^{\mathrm{div}}\Big)\cap L^\infty\Big(0,T;\boldsymbol{H}^{\mathrm{div}}\Big). \tag{8}$$

Its time derivative, however, can generally only be shown to fulfil

$$\partial_t\boldsymbol{u}\in L^{4/d}\Big(0,T;\big(\boldsymbol{V}^{\mathrm{div}}\big)^*\Big). \tag{9}$$

Therefore, (5) represents an assumption for the regularity of $\partial_t\boldsymbol{u}$ both in time (only for $d=3$) and space. The reasons for this simplification are analogous to Remark 2.1.

### 2.2 Essential regularity assumption for computable flows

In addition to the above introduced regularity assumptions which have a direct impact on the weak formulation of the Navier–Stokes problem, it is very common to assume that the solution $\boldsymbol{u}$ to (6) fulfils

$$\nabla \boldsymbol{u} \in L^1(0,T;\boldsymbol{L}^\infty(\Omega)). \tag{10}$$

This short section is aimed at highlighting that incompressible flows which have the essential regularity (10) are relevant both from a theoretical and a practical viewpoint. Let us give a few arguments underlining this statement. At first, (10) guarantees unique solvability of the Navier–Stokes problem; cf. [64, Lemma 2.2]. In fact, (10) ensures that the velocity field $\boldsymbol{u}$ is uniformly Lipschitz continuous on $[0,T]$. As a consequence, the characteristic curves of the dynamical system $\frac{\mathrm{d}}{\mathrm{d}t}\boldsymbol{x}(t) = \boldsymbol{u}(t,\boldsymbol{x}(t))$ remain smooth and never intersect within $[0,T]$; cf. [5]. From a physical point of view, these characteristic curves are the pathlines of the flow; cf. [27, Section 4.3.1]. Lastly, the symmetric part of the velocity gradient $\nabla \boldsymbol{u}$ encodes relevant information about the local structure of a flow; cf. [15, Section 2.5]. In particular, at least in a periodic box and for $\boldsymbol{f} = \boldsymbol{0}$, the smallest scales of an incompressible Navier–Stokes flow behave like $\sqrt{\nu/\|\nabla \boldsymbol{u}\|_{\boldsymbol{L}^\infty}}$; cf. [36].

## 3 Abstract discrete setting and FEM

In this chapter, we attempt to define an abstract discrete setting in which all of the FE methods under consideration can be embedded. To this end, the discrete space-time velocity and pressure spaces are

$$\boldsymbol{V}_h^T = \left\{ \boldsymbol{v}_h \in L^2(0,T;\boldsymbol{V}_h) \colon \partial_t \boldsymbol{v}_h \in L^2(0,T;\boldsymbol{V}_h) \right\} \quad \text{and} \quad Q_h^T = L^2(0,T;Q_h). \tag{11}$$

Contrary to the continuous setting in Section 2, we will not explicitly define the discrete spaces $\boldsymbol{V}_h$ and $Q_h$ at this point. Instead, only general assumptions for the FE pair $\boldsymbol{V}_h/Q_h$ are introduced. Before we begin with the minimal global regularity requirements for the spaces, the following standard decomposition of the domain is introduced.

Let $\mathcal{T}_h$ be a shape-regular FE partition of $\Omega$ without hanging nodes and mesh size $h = \max_{K\in\mathcal{T}_h} h_K$, where $h_K$ denotes the diameter of the particular element $K \in \mathcal{T}_h$. The skeleton $\mathcal{F}_h$ denotes the set of all facets with $\mathcal{F}_K = \{F \in \mathcal{F}_h \colon F \subset \partial K\}$ and $N_\partial = \max_{K\in\mathcal{T}_h} \operatorname{card}(\mathcal{F}_K)$. Moreover, $\mathcal{F}_h = \mathcal{F}_h^i \cup \mathcal{F}_h^\partial$ where $\mathcal{F}_h^i$ is the subset of interior facets and $\mathcal{F}_h^\partial$ collects all boundary facets $F \subset \partial\Omega$. To any $F \in \mathcal{F}_h$ we assign a unit normal vector $\boldsymbol{n}_F$ where, for $F \in \mathcal{F}_h^\partial$, this is the outer unit normal vector $\boldsymbol{n}$. If $F \in \mathcal{F}_h^i$, there are two adjacent elements $K^+$ and $K^-$ sharing the facet $F = \overline{\partial K^+} \cap \overline{\partial K^-}$ and $\boldsymbol{n}_F$ points in an arbitrary but fixed direction. Let $\phi$ be any piecewise smooth (scalar-, vector- or tensor-valued) function with traces from within the interior of $K^\pm$ denoted by $\phi^\pm$, respectively. Then, we define the jump $[\![\cdot]\!]_F$ and average $\{\!\!\{\cdot\}\!\!\}_F$ operator across interior facets $F \in \mathcal{F}_h^i$ by

$$[\![\phi]\!]_F = \phi^+ - \phi^- \quad \text{and} \quad \{\!\!\{\phi\}\!\!\}_F = \frac{1}{2}\left(\phi^+ + \phi^-\right). \tag{12}$$

For boundary facets $F \in \mathcal{F}_h^\partial$ we set $[\![\phi]\!]_F = \{\!\!\{\phi\}\!\!\}_F = \phi$. These operators act componentwise for vector- and tensor-valued functions. Frequently, the subscript indicating the facet is omitted.

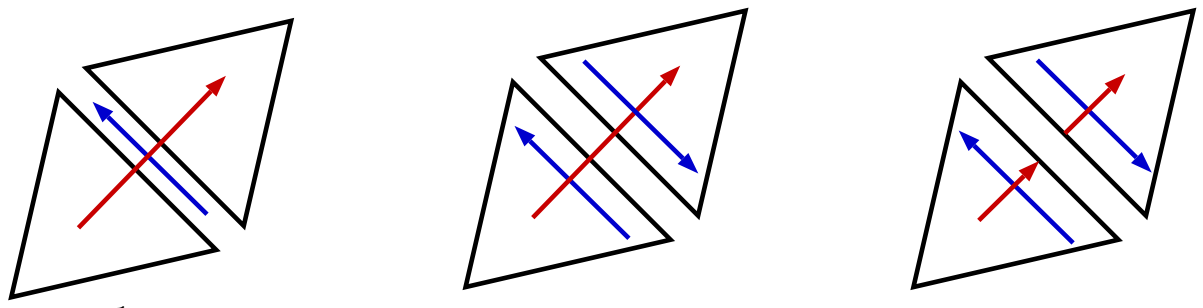

**Fig. 1** Continuity/discontinuity of normal (red) and tangential (blue) components in different methods.

**Assumption A**

$$\boldsymbol{V}_h \subset \boldsymbol{H}(\mathrm{div};\Omega), \quad Q_h \subset L_0^2(\Omega) = Q. \tag{13}$$

Thus, as $\boldsymbol{V} \subset \boldsymbol{H}(\mathrm{div};\Omega)$, our considerations include both $\boldsymbol{H}^1$-conforming and $\boldsymbol{H}(\mathrm{div})$-conforming methods. Fully discontinuous DG-FEM, however, are excluded since $\boldsymbol{L}^2(\Omega) \not\subset \boldsymbol{H}(\mathrm{div};\Omega)$ and at least continuity in normal direction is needed. In Figure 1 a sketch of how the normal and tangential velocity components behave in different methods can be seen.

**Assumption B** The global spaces $\boldsymbol{V}_h$ and $Q_h$ form a discretely inf-sup stable FE pair. That is, there exists $\beta^* > 0$, independent of the mesh size $h$, such that

$$\inf_{q_h \in Q_h \setminus \{0\}} \sup_{\boldsymbol{v}_h \in \boldsymbol{V}_h \setminus \{\boldsymbol{0}\}} \frac{b(\boldsymbol{v}_h, q_h)}{|||\boldsymbol{v}_h|||_e \, \|q_h\|_{L^2}} \geqslant \beta^*. \tag{14}$$

Here, $|||\cdot|||_e$ denotes a suitable energy norm. Due to the $\boldsymbol{H}(\mathrm{div})$-conformity of $\boldsymbol{V}_h$, the pressure-velocity coupling $b(\cdot,\cdot)$ remains the same in the discrete setting. Note that (14) ensures that the space of discretely divergence-free velocities, $\boldsymbol{V}_h^{\mathrm{div}}$, is non-trivial, that is

$$\boldsymbol{V}_h^{\mathrm{div}} = \{\boldsymbol{v}_h \in \boldsymbol{V}_h \colon b(\boldsymbol{v}_h, q_h) = 0, \ \forall q_h \in Q_h\} \neq \{\boldsymbol{0}\}. \tag{15}$$

**Assumption C** The global spaces $\boldsymbol{V}_h$ and $Q_h$ are divergence-conforming, that is

$$\nabla \cdot \boldsymbol{V}_h \subseteq Q_h. \tag{16}$$

If Assumption C holds, the velocity approximation will be exactly divergence-free; cf. [41].

**Remark 3.1** There are several FE pairs which fit into the above introduced framework. The probably most frequently used elements fulfil Assumptions A with $\boldsymbol{V}_h \subset \boldsymbol{H}^1$ and Assumption B. For example, the Taylor–Hood element of order $k$ or the MINI element are well-known; cf. [39] also for different pairs. If Assumption C has to be fulfilled additionally, the Scott–Vogelius element (with certain restrictions on mesh and order) is known. Some other examples are mentioned in [64]. In the context of isogeometric analysis, several $\boldsymbol{H}^1$-conforming and divergence-free FE spaces have been constructed using splines on tensor-product meshes [11,30,31].

Leaving the $\boldsymbol{H}^1$-conforming sector, several classical examples of inf-sup stable $\boldsymbol{H}(\mathrm{div})$-conforming spaces which also fulfil Assumption C can be found in [23,9]. The corresponding methods are discontinuous Galerkin (DG) methods since the tangential components are in general discontinuous across interior facets. Let us specifically mention the family of Raviart–Thomas (RT) elements on simplicial meshes which, for example, have been used in [65]. In this work, however, the family of Brezzi–Douglas–Marini (BDM) elements (applicable on either simplicial or tensor-product meshes) is used in Subsection 5.3. Let us mention that the computational efficiency of $\boldsymbol{H}(\mathrm{div})$-conforming methods can be improved drastically by hybridisation; cf., for example, [46]. Some computational aspects of $\boldsymbol{H}(\mathrm{div})$-conforming FEM are discussed in the Appendix.

### 3.1 Finite element method

The space-semidiscrete (or continuous-in-time) weak formulation of (6) reads as follows:

$$\begin{cases} \text{Find } (\boldsymbol{u}_h, p_h) \in \boldsymbol{V}_h^T \times Q_h^T \text{ with } \boldsymbol{u}_h(0) = \boldsymbol{u}_{0h} \text{ s.t., } \forall (\boldsymbol{v}_h, q_h) \in \boldsymbol{V}_h \times Q_h, & \text{(17a)} \\ (\partial_t \boldsymbol{u}_h, \boldsymbol{v}_h) + \nu a_h(\boldsymbol{u}_h, \boldsymbol{v}_h) + c_h(\boldsymbol{u}_h; \boldsymbol{u}_h, \boldsymbol{v}_h) + b(\boldsymbol{v}_h, p_h) - b(\boldsymbol{u}_h, q_h) = (\boldsymbol{f}, \boldsymbol{v}_h). & \text{(17b)} \end{cases}$$

Note that since the approximation $\boldsymbol{u}_h \in \boldsymbol{V}_h^T$ to (17) does not necessarily have to be $\boldsymbol{H}^1$-conforming, we introduce the broken Sobolev space

$$\boldsymbol{H}^m(\mathcal{T}_h) = \left\{ \boldsymbol{w} \in \boldsymbol{L}^2(\Omega) \colon \boldsymbol{w}\big|_K \in \boldsymbol{H}^m(K),\ \forall K \in \mathcal{T}_h \right\}. \tag{18}$$

Define the broken gradient $\nabla_h \colon \boldsymbol{H}^1(\mathcal{T}_h) \to \boldsymbol{L}^2(\Omega)$ by $(\nabla_h \boldsymbol{w})\big|_K = \nabla\big(\boldsymbol{w}\big|_K\big)$. To be mathematically more precise, the appearance of traces of velocity facet values and normal derivatives thereof dictates that the velocities, at least, belong to $\boldsymbol{H}^{\frac{3}{2}+\varepsilon}(\mathcal{T}_h)$ for some $\varepsilon > 0$; cf. [59, Section 2.1.3].

For the discretisation of the diffusion term, we employ the standard symmetric interior penalty (SIP) form [59,26] (jump penalisation parameter $\sigma > 0$) with an additional grad-div term (parameter $\delta \geqslant 0$):

$$a_h(\boldsymbol{w}, \boldsymbol{v}_h) = \int_\Omega \nabla_h \boldsymbol{w} : \nabla_h \boldsymbol{v}_h \,\mathrm{d}\boldsymbol{x} + \frac{\delta}{\nu} \int_\Omega (\nabla \cdot \boldsymbol{w})(\nabla \cdot \boldsymbol{v}_h) \,\mathrm{d}\boldsymbol{x} \tag{19a}$$

$$- \sum_{F \in \mathcal{F}_h} \oint_F \left[ \{\!\!\{ \nabla \boldsymbol{w} \}\!\!\} \boldsymbol{n}_F \cdot [\![ \boldsymbol{v}_h ]\!] + [\![ \boldsymbol{w} ]\!] \cdot \{\!\!\{ \nabla \boldsymbol{v}_h \}\!\!\} \boldsymbol{n}_F - \frac{\sigma}{h_F} [\![ \boldsymbol{w} ]\!] \cdot [\![ \boldsymbol{v}_h ]\!] \right] \mathrm{d}\boldsymbol{s} \tag{19b}$$

For $\boldsymbol{H}^1$-FEM (globally continuous), the summation over all facets terms $F \in \mathcal{F}_h$ disappears since in this case all jumps vanish. Also, the broken gradient in the first volume term is simply the usual gradient. Whenever the considered FE pair fulfils (16), the discrete velocity is pointwise divergence-free and the grad-div term vanishes. In conjunction with the viscous term $a_h$, the following norms are used:

$$|\!|\!| \boldsymbol{w} |\!|\!|_e^2 = \|\nabla_h \boldsymbol{w}\|_{\boldsymbol{L}^2}^2 + \sum_{F \in \mathcal{F}_h} \frac{\sigma}{h_F} \|[\![ \boldsymbol{w} ]\!]\|_{\boldsymbol{L}^2(F)}^2 \tag{20a}$$

$$|\!|\!| \boldsymbol{w} |\!|\!|_{e,\sharp}^2 = |\!|\!| \boldsymbol{w} |\!|\!|_e^2 + \sum_{K \in \mathcal{T}_h} h_K \|\nabla \boldsymbol{w} \cdot \boldsymbol{n}_K\|_{\boldsymbol{L}^2(\partial K)}^2 \tag{20b}$$

Here, $|\!|\!| \cdot |\!|\!|_e$ denotes the discrete energy norm and the index $\sharp$ indicates a stronger norm. Furthermore, we define the following physically relevant quantities.

**Definition 3.2 (Kinetic and dissipation energies)** The kinetic energy and the kinetic energy dissipation rate of a flow, represented by the velocity $\boldsymbol{w}$, at almost every $t \in (0,T)$ is given, respectively, by

$$\frac{1}{2} \|\boldsymbol{w}(t)\|_{\boldsymbol{L}^2}^2 \quad \text{and} \quad \nu |\!|\!| \boldsymbol{w}(t) |\!|\!|_e^2. \tag{21}$$

For the inertia term, we choose the following convection term [26] for $\boldsymbol{\beta} \in \boldsymbol{L}^\infty \cap \boldsymbol{H}(\mathrm{div};\Omega)$ with $\boldsymbol{\beta} \cdot \boldsymbol{n}\big|_{\partial\Omega} = 0$:

$$c_h(\boldsymbol{\beta};\boldsymbol{w},\boldsymbol{v}_h) = \int_\Omega (\boldsymbol{\beta}\cdot\nabla_h)\boldsymbol{w}\cdot\boldsymbol{v}_h \,\mathrm{d}\boldsymbol{x} + \frac{1}{2}\int_\Omega (\nabla\cdot\boldsymbol{\beta})(\boldsymbol{w}\cdot\boldsymbol{v}_h)\,\mathrm{d}\boldsymbol{x} \tag{22a}$$

$$- \sum_{F\in\mathcal{F}_h^i} \oint_F (\boldsymbol{\beta}\cdot\boldsymbol{n}_F)[\![\boldsymbol{w}]\!]\cdot\{\!\{\boldsymbol{v}_h\}\!\}\,\mathrm{d}\boldsymbol{s} + \sum_{F\in\mathcal{F}_h^i} \oint_F \frac{1}{2}|\boldsymbol{\beta}\cdot\boldsymbol{n}_F|[\![\boldsymbol{w}]\!]\cdot[\![\boldsymbol{v}_h]\!]\,\mathrm{d}\boldsymbol{s} \tag{22b}$$

For $\boldsymbol{H}^1$-conforming FEM the second volume term represents a skew-symmetrisation (other choices are possible; cf. [39]) which vanishes for $\nabla\cdot\boldsymbol{\beta} = 0$. In the general case, the first three terms together are skew-symmetric. $\boldsymbol{H}(\mathrm{div})$-FEM, due to discontinuity in tangential direction, provide the opportunity of including a natural upwind mechanism for stabilising high Reynolds number flows [23]. The corresponding terms are the facet integrals in (22) where the last part is symmetric positive semidefinite. Again, in the globally continuous case, all facet terms vanish. In order to highlight the impact of the upwind term, we introduce the jump seminorm

$$|\boldsymbol{w}|^2_{\boldsymbol{\beta},\mathrm{upw}} = \sum_{F\in\mathcal{F}_h^i} \oint_F \frac{1}{2}|\boldsymbol{\beta}\cdot\boldsymbol{n}_F|\,|[\![\boldsymbol{w}]\!]|^2\,\mathrm{d}\boldsymbol{s}. \tag{23}$$

**Remark 3.3** As can be seen from (19) and (22), an exactly divergence-free and $\boldsymbol{H}^1$-conforming method leads to a scheme which, in terms of multilinear forms, is identical to the continuous one in (6b). In this sense, divergence-free $\boldsymbol{H}^1$-FEM represent, at least from a theoretical point of view, the most simplified available FE methods. Hence, the numerical analysis for this class of methods is also the most concise and compact.

## 3.2 Energy estimate and well-posedness

Let us summarise the most important discrete coercivity properties; cf. [59, 26].

**Lemma 3.4 (Discrete coercivity of $a_h$ and $c_h$)** *Assume that $\sigma > 0$ is sufficiently large. Then, the bilinear form $a_h$ is coercive on $\boldsymbol{V}_h$ w.r.t. the energy norm $|\!|\!|\cdot|\!|\!|_e$. Moreover, the grad-div term allows for an additional control over the divergence of the discrete velocity. The convective form $c_h$ is coercive on $\boldsymbol{V}_h$ w.r.t. the upwind seminorm $|\cdot|_{\mathrm{upw}}$. That is, there exists $C_\sigma > 0$, independent of $h$, such that, for all $\boldsymbol{v}_h \in \boldsymbol{V}_h$,*

$$a_h(\boldsymbol{v}_h,\boldsymbol{v}_h) \geqslant C_\sigma |\!|\!|\boldsymbol{v}_h|\!|\!|_e^2 + \frac{\delta}{\nu}\|\nabla\cdot\boldsymbol{v}_h\|^2_{\boldsymbol{L}^2} \quad \text{and} \quad c_h(\boldsymbol{\beta};\boldsymbol{v}_h,\boldsymbol{v}_h) = |\boldsymbol{v}_h|^2_{\boldsymbol{\beta},\mathrm{upw}}. \tag{24}$$

As a consequence, and after applying standard arguments, we obtain the following corollary.

**Corollary 3.5 (Well-posedness and velocity energy estimate)** *Let $\boldsymbol{f} \in L^1(\boldsymbol{L}^2)$ and $\boldsymbol{u}_{0h} \in \boldsymbol{L}^2$. Then, there exists a solution $\boldsymbol{u}_h \in \boldsymbol{V}_h^T$ to (17) with*

$$\frac{1}{2}\|\boldsymbol{u}_h(T)\|^2_{\boldsymbol{L}^2} + \int_0^T \left[\nu C_\sigma |\!|\!|\boldsymbol{u}_h|\!|\!|_e^2 + \delta\|\nabla\cdot\boldsymbol{u}_h\|^2_{\boldsymbol{L}^2} + |\boldsymbol{u}_h|^2_{\boldsymbol{u}_h,\mathrm{upw}}\right]\mathrm{d}\tau \tag{25a}$$

$$\leqslant \|\boldsymbol{u}_{0h}\|^2_{\boldsymbol{L}^2} + \frac{3}{2}\|\boldsymbol{f}\|^2_{L^1(\boldsymbol{L}^2)}. \tag{25b}$$

*Provided $\boldsymbol{f}$ is even Lipschitz in time, the solution $\boldsymbol{u}_h$ is unique.*

## 4 Classical inf-sup stable $H^1$-conforming FEM

In this section we briefly want to discuss non-divergence-free $\boldsymbol{H}^1$-conforming methods.

### 4.1 Classical $\boldsymbol{H}^1$-conforming mixed FEM are not *Re*-semi-robust

The application of the Gronwall lemma to continuous-in-time estimates of both kinetic and dissipation energies for the Galerkin-FEM leads to an exponential factor on the right-hand side which may depend in an unfavourable way on the length of the time interval $(0,T)$, norms of the solution, and on inverse powers of the viscosity.

More precisely, for the skew-symmetric form (22) of the convective term and only assuming $\nabla \boldsymbol{u} \in L^4(\boldsymbol{L}^2)$, the argument of the exponential on the right-hand side is of the form $C\nu^{-3} \|\nabla \boldsymbol{u}\|^4_{L^4(\boldsymbol{L}^2)}$; cf. [39, Theorem 7.35]. Following [39, Remark 7.39], under the assumptions $\nabla \boldsymbol{u} \in L^1(\boldsymbol{L}^\infty)$ and $\boldsymbol{u} \in L^2(\boldsymbol{L}^\infty)$, one can improve the argument of the exponential to

$$\frac{1}{2} \|\nabla \boldsymbol{u}\|_{L^1(\boldsymbol{L}^\infty)} + \frac{4}{\nu} \|\boldsymbol{u}\|^2_{L^2(\boldsymbol{L}^\infty)}. \tag{26}$$

Nevertheless, estimates with such strong exponential growth are useless in practice. Please note that even a rough error estimate using the triangle inequality together with the stability estimates for both discrete and continuous solution, see Corollary 3.5, provides asymptotically much better *Re*-semi-robust bounds as opposed to exponential growth depending on $\nu^{-1}$. We refer to Subsection 5.3 where it can be observed that the error from a classical Taylor–Hood Galerkin computation actually shows such an unfavourable exponential growth.

### 4.2 Improvements with grad-div stabilisation

To the best of our knowledge, it has been observed first in [51] that for inf-sup stable FE pairs, the combination of the Galerkin-FEM with grad-div stabilisation can avoid entirely the explicit dependence of the Gronwall factor on $\nu^{-1}$. For $\boldsymbol{H}^1$-conforming inf-sup stable FEM, this provides *Re*-semi-robust estimates; see the results [3] which were improved in [25]. In particular, the dependence of the argument of the Gronwall factor (26) can be replaced by

$$T + C_1 \|\nabla \boldsymbol{u}\|_{L^1(\boldsymbol{L}^\infty)} + C_2 \frac{h^2}{\delta} \|\boldsymbol{u}\|^2_{L^1(\boldsymbol{W}^{1,\infty})}, \tag{27}$$

where $\delta$ is the grad-div parameter. The work [33] deepens the results; in particular for optimal pressure estimates. For the argument of the Gronwall factor, they obtain a Gronwall argument similar to (27) where the third summand is replaced by $\frac{1}{2\delta} \|\boldsymbol{u}\|^2_{L^1(\boldsymbol{H}^2)}$.

Our numerical experience shows that grad-div stabilisation can improve the results for classical inf-sup stable $\boldsymbol{H}^1$-conforming Galerkin-FEM. For example, for the numerical simulation of a problem with standing vortices in Subsection 5.3, we show that sometimes grad-div stabilisation dramatically improves the behaviour of the Gronwall factor. However, recent research has clarified that these effects are not really a stabilisation issue, but are related to some kind of consistency error, whenever the fundamental invariance property (replacing the source term $\boldsymbol{f}$ by $\boldsymbol{f}+\nabla\psi$ changes the solution $(\boldsymbol{u},p)$ to $(\boldsymbol{u},p+\psi)$) is violated on the discrete level. Therefore, grad-div stabilisation simply reduces the divergence error of discrete velocity solutions, which involves the potential danger of the classical Poisson locking phenomenon [38]. Also, further aspects of the numerical analysis of grad-div stabilisation of such FEM can be found in the work [14]. These considerations made us choose inf-sup stable, exactly divergence-free mixed FEM for this article; see Section 5.

## 5 Divergence-free $H^1$- and $H$(div)-FEM

Under Assumption C, the following Galerkin orthogonality property in $\boldsymbol{V}_h^{\mathrm{div}}$ can be stated without any contributions from the pressure. The most important ingredient is the consistency of both SIP formulation of the viscous term and upwind formulation of the convective term [59,26].

**Corollary 5.1 (Galerkin orthogonality)** *Let $\boldsymbol{u}_h \in \boldsymbol{V}_h^T$ solve (17). Assume that the solution $\boldsymbol{u} \in \boldsymbol{V}^T$ of (6) satisfies the regularity condition $\boldsymbol{u} \in L^2\left(\boldsymbol{H}^{\frac{3}{2}+\varepsilon}(\mathcal{T}_h)\right)$ for $\varepsilon > 0$. Then, for a.e. $t \in (0,T)$ and for all $\boldsymbol{v}_h \in \boldsymbol{V}_h^{\mathrm{div}}$,*

$$(\partial_t[\boldsymbol{u}-\boldsymbol{u}_h],\boldsymbol{v}_h)+\nu a_h(\boldsymbol{u}-\boldsymbol{u}_h,\boldsymbol{v}_h)+c_h(\boldsymbol{u};\boldsymbol{u},\boldsymbol{v}_h)-c_h(\boldsymbol{u}_h;\boldsymbol{u}_h,\boldsymbol{v}_h)=0. \tag{28}$$

### 5.1 Stationary Stokes projection

In this section we want to consider the coupling of pressure and viscous effects only. With a sufficiently smooth forcing term $\boldsymbol{g}$, the well-known continuous weak formulation of the stationary Stokes problem reads

$$\begin{cases} \text{Find } (\boldsymbol{u}_s,p_s) \in \boldsymbol{V}\times Q \text{ s.t., } \forall (\boldsymbol{v},q) \in \boldsymbol{V}\times Q, & \text{(29a)} \\ \nu a(\boldsymbol{u}_s,\boldsymbol{v})+b(\boldsymbol{v},p_s)-b(\boldsymbol{u}_s,q)=(\boldsymbol{g},\boldsymbol{v}). & \text{(29b)} \end{cases}$$

In order to obtain optimal $\boldsymbol{L}^2$-estimates for the velocity, we make the following assumption which is called 'elliptic regularity', 'Cattabriga's regularity' or 'smoothing property'.

**Assumption D** Assume that $\Omega$ is either a convex polygon for $d=2$ or of class $\mathcal{C}^{1,1}$ for $d\in\{2,3\}$. Then, for all $\boldsymbol{g}\in\boldsymbol{L}^2$, the solution $(\boldsymbol{u}_s,p_s)\in\boldsymbol{V}\times Q$ of (29) additionally fulfils the regularity property $(\boldsymbol{u}_s,p_s)\in\boldsymbol{H}^2\times H^1$ and the energy estimate $\sqrt{\nu}\,\|\boldsymbol{u}_s\|_{\boldsymbol{H}^2}+\|p_s\|_{H^1}\leqslant C\,\|\boldsymbol{g}\|_{\boldsymbol{L}^2}$; cf. [10, Theorem IV.5.8].

Note that the following definition is stated directly in $\boldsymbol{V}_h^{\mathrm{div}}$ because this suffices for our considerations.

**Definition 5.2 (Stationary Stokes projection)** Let $\boldsymbol{w}\in\boldsymbol{H}^{\frac{3}{2}+\varepsilon}(\mathcal{T}_h)$ for $\varepsilon>0$ fulfil $\nabla\cdot\boldsymbol{w}=0$ pointwise. Then, we define the stationary Stokes projection $\boldsymbol{\pi}_s\boldsymbol{w}\in\boldsymbol{V}_h^{\mathrm{div}}$ of $\boldsymbol{w}$ to be the unique FE solution to the problem

$$a_h(\boldsymbol{\pi}_s\boldsymbol{w},\boldsymbol{v}_h)=a_h(\boldsymbol{w},\boldsymbol{v}_h),\quad \forall\boldsymbol{v}_h\in\boldsymbol{V}_h^{\mathrm{div}}. \tag{30}$$

As a consequence, the approximation properties of the projection operator $\boldsymbol{\pi}_s$ can be derived from error estimates for the stationary Stokes problem. The following theorem holds true; cf. [64,65].

**Theorem 5.3 (Stokes projection error estimate)** *Let $\boldsymbol{\pi}_s \boldsymbol{w}$ be the Stokes projection of $\boldsymbol{w}$ with $\nabla\cdot\boldsymbol{w} = 0$ and Assumption D be fulfilled. Then, provided $\boldsymbol{w} \in \boldsymbol{H}^r(\Omega)$ with $r > 3/2$ and $r_{\boldsymbol{u}} = \min\{r, k+1\}$,*

$$\|\boldsymbol{w}-\boldsymbol{\pi}_s\boldsymbol{w}\|_{\boldsymbol{L}^2} + h\|\!|\boldsymbol{w}-\boldsymbol{\pi}_s\boldsymbol{w}\|\!|_{e,\sharp} \leqslant Ch \inf_{\boldsymbol{v}_h \in \boldsymbol{V}_h^{\mathrm{div}}} \|\!|\boldsymbol{w}-\boldsymbol{v}_h\|\!|_{e,\sharp} \leqslant Ch^{r_{\boldsymbol{u}}} |\boldsymbol{w}|_{\boldsymbol{H}^{r_{\boldsymbol{u}}}}. \tag{31}$$

**Assumption E** In the setting of Theorem 5.3, depending on which method is used, we assume that

$$\boldsymbol{H}^1: \quad \|\nabla_h \boldsymbol{\pi}_s \boldsymbol{w}\|_{\boldsymbol{L}^\infty} \leqslant C\|\nabla_h \boldsymbol{w}\|_{\boldsymbol{L}^\infty}, \tag{32a}$$

$$\boldsymbol{H}(\mathrm{div}): \quad \|\boldsymbol{w}-\boldsymbol{\pi}_s\boldsymbol{w}\|_{\boldsymbol{L}^\infty} + h\|\nabla_h\boldsymbol{\pi}_s\boldsymbol{w}\|_{\boldsymbol{L}^\infty} \leqslant Ch\|\nabla_h\boldsymbol{w}\|_{\boldsymbol{L}^\infty}. \tag{32b}$$

**Remark 5.4** In the $\boldsymbol{H}^1$-conforming context, an analogue to (32a) has been shown in [34] in the context of non-divergence-free methods which involves also the pressure. The analysis in [34] simplifies for divergence-free $\boldsymbol{H}^1$-conforming methods, thereby leading to (32a). The validity of (32b) is an open problem although, in principle, similar techniques as in [34] seem to be applicable. We are not aware of any literature where $\boldsymbol{L}^\infty$ estimates for the $\boldsymbol{H}(\mathrm{div})$-conforming Stokes projection have been discussed. Note that in [35] the assumption (32b) is circumvented by assuming a similar estimate for an $\boldsymbol{H}(\mathrm{div})$-conforming interpolation.

## 5.2 Pressure- and *Re*-semi-robust error estimates

In this section, additionally to pressure and viscous effects, the dynamics of the Navier–Stokes problem are investigated; this means the evolutionary and inertia term. We use the Stokes projection to introduce the error splitting

$$\boldsymbol{u}-\boldsymbol{u}_h = [\boldsymbol{u}-\boldsymbol{\pi}_s\boldsymbol{u}] - [\boldsymbol{u}_h - \boldsymbol{\pi}_s\boldsymbol{u}] = \boldsymbol{\eta} - \boldsymbol{e}_h. \tag{33}$$

For the $\boldsymbol{H}(\mathrm{div})$-conforming methods we additionally need to be able to bound facet norms by volume norms:

**Assumption F** The velocity space $\boldsymbol{V}_h$ satisfies the discrete trace inequality [26, Remark 1.47]

$$\forall \boldsymbol{v}_h \in \boldsymbol{V}_h: \quad \|\boldsymbol{v}_h\|_{\boldsymbol{L}^2(\partial K)} \leqslant C_{\mathrm{tr}} h_K^{-1/2} \|\boldsymbol{v}_h\|_{\boldsymbol{L}^2(K)}, \quad \forall K \in \mathcal{T}_h. \tag{34}$$

**Lemma 5.5 (Difference of convective terms)** *Assume that $\boldsymbol{u} \in L^1(\boldsymbol{W}^{1,\infty})$. Then, for all finite $\varepsilon_i > 0$, $i \in \{1,2,3,4\}$, we obtain*

$$c_h(\boldsymbol{u};\boldsymbol{u},\boldsymbol{e}_h) - c_h(\boldsymbol{u}_h;\boldsymbol{u}_h,\boldsymbol{e}_h) \leqslant -|\boldsymbol{e}_h|^2_{\boldsymbol{u}_h,\mathrm{upw}} \tag{35a}$$

$$+\varepsilon_1^{-1}\|\boldsymbol{u}\|_{\boldsymbol{L}^\infty}\|\nabla_h\boldsymbol{\eta}\|^2_{\boldsymbol{L}^2} + \left[\varepsilon_2^{-1} + C\left(\varepsilon_3^{-1}+\varepsilon_4^{-1}\right)h^{-1}\right]\|\nabla\boldsymbol{u}\|_{\boldsymbol{L}^\infty}\|\boldsymbol{\eta}\|^2_{\boldsymbol{L}^2} \tag{35b}$$

$$+\left[\varepsilon_1\|\boldsymbol{u}\|_{\boldsymbol{L}^\infty} + \left(C+\varepsilon_2+C(\varepsilon_3+\varepsilon_4)h^{-1}\right)\|\nabla\boldsymbol{u}\|_{\boldsymbol{L}^\infty}\right]\|\boldsymbol{e}_h\|^2_{\boldsymbol{L}^2}. \tag{35c}$$

**Proof:** We basically follow the ideas from [35]. At first, insert the definition of $c_h$, use $[\![\boldsymbol{u}]\!]_F = \boldsymbol{0}$ for all facets $F \in \mathcal{F}_h^i$ and reorder:

$$c_h(\boldsymbol{u};\boldsymbol{u},\boldsymbol{e}_h) - c_h(\boldsymbol{u}_h;\boldsymbol{u}_h,\boldsymbol{e}_h) = \int_\Omega \left[(\boldsymbol{u}\cdot\nabla_h)\boldsymbol{u}\cdot\boldsymbol{e}_h - (\boldsymbol{u}_h\cdot\nabla_h)\boldsymbol{u}_h\cdot\boldsymbol{e}_h\right]\mathrm{d}\boldsymbol{x} \tag{36a}$$

$$-\sum_{F\in\mathcal{F}_h^i}\oint_F (\boldsymbol{u}_h\cdot\boldsymbol{n}_F)[\![\boldsymbol{u}-\boldsymbol{u}_h]\!]\cdot\{\!\!\{\boldsymbol{e}_h\}\!\!\}\,\mathrm{d}\boldsymbol{s} + \sum_{F\in\mathcal{F}_h^i}\oint_F \frac{1}{2}|\boldsymbol{u}_h\cdot\boldsymbol{n}_F|[\![\boldsymbol{u}-\boldsymbol{u}_h]\!]\cdot[\![\boldsymbol{e}_h]\!]\,\mathrm{d}\boldsymbol{s} \tag{36b}$$

$$=\mathfrak{T}_1+\mathfrak{T}_2+\mathfrak{T}_3 \tag{36c}$$

Note that in the $\boldsymbol{H}^1$-conforming case, $\mathfrak{T}_2 = \mathfrak{T}_3 = 0$. For the volume term $\mathfrak{T}_1$, we subtract and add $(\boldsymbol{u}\cdot\nabla_h \boldsymbol{\pi}_s\boldsymbol{u}, \boldsymbol{e}_h)_\Omega$, replace $\boldsymbol{u}_h = \boldsymbol{\pi}_s\boldsymbol{u} + \boldsymbol{e}_h$ and use triangle, Hölder's and Young's $(\varepsilon_1, \varepsilon_2 > 0)$ inequality:

$$
\begin{aligned}
\mathfrak{T}_1 + (\boldsymbol{u}_h\cdot\nabla_h\boldsymbol{e}_h, \boldsymbol{e}_h)_\Omega &= (\boldsymbol{u}\cdot\nabla_h\boldsymbol{\eta}, \boldsymbol{e}_h)_\Omega + ([\boldsymbol{u}-\boldsymbol{u}_h]\cdot\nabla_h\boldsymbol{\pi}_s\boldsymbol{u}, \boldsymbol{e}_h)_\Omega && \text{(37a)}\\
&\leqslant \|\boldsymbol{u}\|_{\boldsymbol{L}^\infty}\|\nabla_h\boldsymbol{\eta}\|_{\boldsymbol{L}^2}\|\boldsymbol{e}_h\|_{\boldsymbol{L}^2} + \|\boldsymbol{\eta}-\boldsymbol{e}_h\|_{\boldsymbol{L}^2}\|\nabla_h\boldsymbol{\pi}_s\boldsymbol{u}\|_{\boldsymbol{L}^\infty}\|\boldsymbol{e}_h\|_{\boldsymbol{L}^2} && \text{(37b)}\\
&\leqslant \varepsilon_1^{-1}\|\boldsymbol{u}\|_{\boldsymbol{L}^\infty}\|\nabla_h\boldsymbol{\eta}\|_{\boldsymbol{L}^2}^2 + \varepsilon_2^{-1}\|\nabla_h\boldsymbol{\pi}_s\boldsymbol{u}\|_{\boldsymbol{L}^\infty}\|\boldsymbol{\eta}\|_{\boldsymbol{L}^2}^2 && \text{(37c)}\\
&+ [\varepsilon_1\|\boldsymbol{u}\|_{\boldsymbol{L}^\infty} + (1+\varepsilon_2)\|\nabla_h\boldsymbol{\pi}_s\boldsymbol{u}\|_{\boldsymbol{L}^\infty}]\|\boldsymbol{e}_h\|_{\boldsymbol{L}^2}^2 && \text{(37d)}
\end{aligned}
$$

For $\boldsymbol{H}^1$-conforming methods, $(\boldsymbol{u}_h\cdot\nabla_h\boldsymbol{e}_h, \boldsymbol{e}_h)_\Omega = (\boldsymbol{u}_h\cdot\nabla\boldsymbol{e}_h, \boldsymbol{e}_h)_\Omega = 0$ and the proof is already complete at this point. For $\boldsymbol{H}(\mathrm{div})$-FEM, this term is balanced by the facet terms. In fact, for these facet terms inserting the error splitting leads to

$$
\begin{aligned}
\mathfrak{T}_2 &= -\sum_{F\in\mathcal{F}_h^i}\oint_F (\boldsymbol{u}_h\cdot\boldsymbol{n}_F)[\![\boldsymbol{\eta}]\!]\cdot\{\!\!\{\boldsymbol{e}_h\}\!\!\}\,\mathrm{d}\boldsymbol{s} + \sum_{F\in\mathcal{F}_h^i}\oint_F (\boldsymbol{u}_h\cdot\boldsymbol{n}_F)[\![\boldsymbol{e}_h]\!]\cdot\{\!\!\{\boldsymbol{e}_h\}\!\!\}\,\mathrm{d}\boldsymbol{s} && \text{(38a)}\\
&= \mathfrak{T}_{2,1} + \mathfrak{T}_{2,2} && \text{(38b)}\\
\mathfrak{T}_3 &= \sum_{F\in\mathcal{F}_h^i}\oint_F \frac{1}{2}|\boldsymbol{u}_h\cdot\boldsymbol{n}_F|[\![\boldsymbol{\eta}]\!]\cdot[\![\boldsymbol{e}_h]\!]\,\mathrm{d}\boldsymbol{s} - \sum_{F\in\mathcal{F}_h^i}\oint_F \frac{1}{2}|\boldsymbol{u}_h\cdot\boldsymbol{n}_F|[\![\boldsymbol{e}_h]\!]\cdot[\![\boldsymbol{e}_h]\!]\,\mathrm{d}\boldsymbol{s} && \text{(38c)}\\
&= \mathfrak{T}_{3,1} - |\boldsymbol{e}_h|_{\boldsymbol{u}_h,\mathrm{upw}}^2 && \text{(38d)}
\end{aligned}
$$

Here, due to the discrete coercivity of $c_h$ (Lemma 3.4), we can conclude that $\mathfrak{T}_{2,2} = (\boldsymbol{u}_h\cdot\nabla_h\boldsymbol{e}_h, \boldsymbol{e}_h)_\Omega$ and thus, in the end, the term cancels out with its corresponding part from the volume term $\mathfrak{T}_1$. For the remaining two facet terms, apply Hölder's inequality after again inserting the relation $\boldsymbol{u}_h = \boldsymbol{e}_h + \boldsymbol{\pi}_s\boldsymbol{u}$:

$$
\begin{aligned}
|\mathfrak{T}_{2,1}| &\leqslant \sum_{F\in\mathcal{F}_h^i}\oint_F \left|(\boldsymbol{e}_h\cdot\boldsymbol{n}_F)[\![\boldsymbol{\eta}]\!]\cdot\{\!\!\{\boldsymbol{e}_h\}\!\!\}\right|\mathrm{d}\boldsymbol{s} + \sum_{F\in\mathcal{F}_h^i}\oint_F \left|(\boldsymbol{\pi}_s\boldsymbol{u}\cdot\boldsymbol{n}_F)[\![\boldsymbol{\eta}]\!]\cdot\{\!\!\{\boldsymbol{e}_h\}\!\!\}\right|\mathrm{d}\boldsymbol{s} && \text{(39a)}\\
&\leqslant \|\boldsymbol{\eta}\|_{\boldsymbol{L}^\infty}\sum_{F\in\mathcal{F}_h^i}\left\|\{\!\!\{\boldsymbol{e}_h\}\!\!\}\right\|_{\boldsymbol{L}^2(F)}^2 + \|\boldsymbol{\pi}_s\boldsymbol{u}\|_{\boldsymbol{L}^\infty}\sum_{F\in\mathcal{F}_h^i}\|[\![\boldsymbol{\eta}]\!]\|_{\boldsymbol{L}^2(F)}\left\|\{\!\!\{\boldsymbol{e}_h\}\!\!\}\right\|_{\boldsymbol{L}^2(F)} && \text{(39b)}\\
&= \mathfrak{T}_{2,1,1} + \mathfrak{T}_{2,1,2} && \text{(39c)}
\end{aligned}
$$

Using the bound $\frac{1}{2}(a+b)^2 \leqslant (a^2+b^2)$ for $a,b\in\mathbb{R}$ and the discrete trace inequality (Assumption F) we observe that

$$
\begin{aligned}
\sum_{F\in\mathcal{F}_h^i}\oint_F \left|\{\!\!\{\boldsymbol{e}_h\}\!\!\}\right|^2\mathrm{d}\boldsymbol{s} &\leqslant \sum_{F\in\mathcal{F}_h^i}\left[\|\boldsymbol{e}_h^+\|_{\boldsymbol{L}^2(F)}^2 + \|\boldsymbol{e}_h^-\|_{\boldsymbol{L}^2(F)}^2\right] && \text{(40a)}\\
&\leqslant \sum_{K\in\mathcal{T}_h}\|\boldsymbol{e}_h\|_{\boldsymbol{L}^2(\partial K)}^2 \leqslant C_{\mathrm{tr}}^2 h^{-1}\|\boldsymbol{e}_h\|_{\boldsymbol{L}^2}^2. && \text{(40b)}
\end{aligned}
$$

The same estimate can be obtained when the average is replaced by the jump over facets. Together with the $\boldsymbol{L}^\infty$ approximation properties of $\boldsymbol{\pi}_s\boldsymbol{u}$ (Assumption E), this results in

$$
|\mathfrak{T}_{2,1,1}| \leqslant C\|\nabla\boldsymbol{u}\|_{\boldsymbol{L}^\infty}\|\boldsymbol{e}_h\|_{\boldsymbol{L}^2}^2. \tag{41}
$$

Similarly, with Young's inequality ($\varepsilon_3 > 0$),

$$|\mathfrak{T}_{2,1,2}| \leqslant \|\boldsymbol{\pi}_s\boldsymbol{u}\|_{\boldsymbol{L}^\infty} \left(\sum_{F\in\mathcal{F}_h^i} \|[\![\boldsymbol{\eta}]\!]\|^2_{\boldsymbol{L}^2(F)}\right)^{1/2} \left(\sum_{F\in\mathcal{F}_h^i} \|\{\!\{\boldsymbol{e}_h\}\!\}\|^2_{\boldsymbol{L}^2(F)}\right)^{1/2} \tag{42a}$$

$$\leqslant C\,\|\boldsymbol{\pi}_s\boldsymbol{u}\|_{\boldsymbol{L}^\infty}\, h^{-1/2}\,\|\boldsymbol{\eta}\|_{\boldsymbol{L}^2}\, h^{-1/2}\,\|\boldsymbol{e}_h\|_{\boldsymbol{L}^2} \tag{42b}$$

$$\leqslant C\varepsilon_3^{-1}\,\|\boldsymbol{\pi}_s\boldsymbol{u}\|_{\boldsymbol{L}^\infty} h^{-1}\,\|\boldsymbol{\eta}\|^2_{\boldsymbol{L}^2} + C\varepsilon_3\,\|\boldsymbol{\pi}_s\boldsymbol{u}\|_{\boldsymbol{L}^\infty} h^{-1}\,\|\boldsymbol{e}_h\|^2_{\boldsymbol{L}^2}. \tag{42c}$$

The estimate of the upwind term $\mathfrak{T}_{3,1}$ is completely analogous after using the triangle inequality in the form $|\boldsymbol{u}_h\cdot\boldsymbol{n}_F| \leqslant |\boldsymbol{e}_h\cdot\boldsymbol{n}_F| + |\boldsymbol{\pi}_s\boldsymbol{u}\cdot\boldsymbol{n}_F|$. With $\varepsilon_4 > 0$, we obtain

$$|\mathfrak{T}_{3,1}| \leqslant C|\boldsymbol{u}|_{\boldsymbol{W}^{1,\infty}}\,\|\boldsymbol{e}_h\|^2_{\boldsymbol{L}^2} + C\varepsilon_4^{-1}\,\|\boldsymbol{\pi}_s\boldsymbol{u}\|_{\boldsymbol{L}^\infty} h^{-1}\,\|\boldsymbol{\eta}\|^2_{\boldsymbol{L}^2} + C\varepsilon_4\,\|\boldsymbol{\pi}_s\boldsymbol{u}\|_{\boldsymbol{L}^\infty} h^{-1}\,\|\boldsymbol{e}_h\|^2_{\boldsymbol{L}^2}. \tag{43}$$

Finally, Assumption E implies stability of the Stokes projection in the form $\|\boldsymbol{\pi}_s\boldsymbol{u}\|_{\boldsymbol{L}^\infty} + \|\nabla_h\boldsymbol{\pi}_s\boldsymbol{u}\|_{\boldsymbol{L}^\infty} \leqslant C|\boldsymbol{u}|_{\boldsymbol{W}^{1,\infty}}$. Combining the above estimates concludes the proof. ■

**Theorem 5.6 (Velocity discretisation error estimate)** *Let $\boldsymbol{u}\in\boldsymbol{V}^T$ solve (6) and $\boldsymbol{u}_h\in\boldsymbol{V}_h^T$ solve (17). If additionally $\boldsymbol{u}\in L^2\left(\boldsymbol{H}^{\frac{3}{2}+\varepsilon}(\mathcal{T}_h)\right)$ for $\varepsilon>0$, $\boldsymbol{u}\in L^1\left(\boldsymbol{W}^{1,\infty}\right)$ and $\boldsymbol{u}_h(0)=\boldsymbol{\pi}_s\boldsymbol{u}_0$, we obtain the following error estimate:*

$$\frac{1}{2}\,\|\boldsymbol{e}_h\|^2_{L^\infty(\boldsymbol{L}^2)} + \int_0^T \left[\nu C_\sigma |\!|\!|\boldsymbol{e}_h|\!|\!|^2_e + |\boldsymbol{e}_h|^2_{\boldsymbol{u}_h,\mathrm{upw}}\right]\mathrm{d}\tau \tag{44a}$$

$$\leqslant e^{G_{\boldsymbol{u}}(T)} \int_0^T \left[\|\partial_t\boldsymbol{\eta}\|^2_{\boldsymbol{L}^2} + \|\boldsymbol{u}\|_{\boldsymbol{L}^\infty}\,\|\nabla_h\boldsymbol{\eta}\|^2_{\boldsymbol{L}^2} + \left(1+Ch^{-2}\right)\|\nabla\boldsymbol{u}\|_{\boldsymbol{L}^\infty}\,\|\boldsymbol{\eta}\|^2_{\boldsymbol{L}^2}\right]\mathrm{d}\tau \tag{44b}$$

*Here, the Gronwall constant is given by*

$$G_{\boldsymbol{u}}(T) = T + \|\boldsymbol{u}\|_{L^1(0,T;\boldsymbol{L}^\infty(\Omega))} + C\,\|\nabla\boldsymbol{u}\|_{L^1(0,T;\boldsymbol{L}^\infty(\Omega))}. \tag{45}$$

**Proof :** Corollary 5.1 with $\boldsymbol{v}_h = \boldsymbol{e}_h(t)\in\boldsymbol{V}_h^{\mathrm{div}}$ and the error splitting (33) yields

$$(\partial_t\boldsymbol{e}_h,\boldsymbol{e}_h) + \nu a_h(\boldsymbol{e}_h,\boldsymbol{e}_h) = (\partial_t\boldsymbol{\eta},\boldsymbol{e}_h) + \nu a_h(\boldsymbol{\eta},\boldsymbol{e}_h) + c_h(\boldsymbol{u};\boldsymbol{u},\boldsymbol{e}_h) - c_h(\boldsymbol{u}_h;\boldsymbol{u}_h,\boldsymbol{e}_h). \tag{46}$$

We use $(\partial_t\boldsymbol{e}_h,\boldsymbol{e}_h) = \frac{1}{2}\frac{\mathrm{d}}{\mathrm{d}t}\,\|\boldsymbol{e}_h\|^2_{\boldsymbol{L}^2}$ and discrete coercivity of $a_h$ (Lemma 3.4) on the left-hand side (note that $\nabla\cdot\boldsymbol{e}_h = 0$). On the right-hand side, apply Cauchy–Schwarz plus Young ($\varepsilon_5>0$) and use Definition 5.2. Then, we obtain

$$\frac{1}{2}\frac{\mathrm{d}}{\mathrm{d}t}\,\|\boldsymbol{e}_h\|^2_{\boldsymbol{L}^2} + \nu C_\sigma |\!|\!|\boldsymbol{e}_h|\!|\!|^2_e \leqslant \varepsilon_5^{-1}\,\|\partial_t\boldsymbol{\eta}\|^2_{\boldsymbol{L}^2} + \varepsilon_5\,\|\boldsymbol{e}_h\|^2_{\boldsymbol{L}^2} + c_h(\boldsymbol{u};\boldsymbol{u},\boldsymbol{e}_h) - c_h(\boldsymbol{u}_h;\boldsymbol{u}_h,\boldsymbol{e}_h). \tag{47}$$

The application of Lemma 5.5 results in

$$\frac{1}{2}\frac{\mathrm{d}}{\mathrm{d}t}\,\|\boldsymbol{e}_h\|^2_{\boldsymbol{L}^2} + \nu C_\sigma |\!|\!|\boldsymbol{e}_h|\!|\!|^2_e + |\boldsymbol{e}_h|^2_{\boldsymbol{u}_h,\mathrm{upw}} \leqslant \varepsilon_5^{-1}\,\|\partial_t\boldsymbol{\eta}\|^2_{\boldsymbol{L}^2} + \varepsilon_1^{-1}\,\|\boldsymbol{u}\|_{\boldsymbol{L}^\infty}\,\|\nabla_h\boldsymbol{\eta}\|^2_{\boldsymbol{L}^2} \tag{48a}$$

$$+ \left[\varepsilon_2^{-1} + C\left(\varepsilon_3^{-1}+\varepsilon_4^{-1}\right)h^{-1}\right]\|\nabla\boldsymbol{u}\|_{\boldsymbol{L}^\infty}\,\|\boldsymbol{\eta}\|^2_{\boldsymbol{L}^2} \tag{48b}$$

$$+ \left[\varepsilon_5 + \varepsilon_1\,\|\boldsymbol{u}\|_{\boldsymbol{L}^\infty} + \left(C+\varepsilon_2+C(\varepsilon_3+\varepsilon_4)h^{-1}\right)\|\nabla\boldsymbol{u}\|_{\boldsymbol{L}^\infty}\right]\|\boldsymbol{e}_h\|^2_{\boldsymbol{L}^2}. \tag{48c}$$

The next step is choosing the $\varepsilon_i$. Note that in this step, numerous different error estimates can be obtained. In the end, everything multiplying $\|\boldsymbol{e}_h\|^2_{\boldsymbol{L}^2}$ will enter the Gronwall exponent and since we do not want to have negative exponents of $h$ there, choosing $\varepsilon_3$ and $\varepsilon_4$ such that $\varepsilon_3 = \varepsilon_4 = \mathcal{O}(h)$ is a valid strategy. For the remaining variables, we simply set $\varepsilon_1 = \varepsilon_2 = \varepsilon_5 = 1$. This results in

$$\frac{1}{2}\frac{\mathrm{d}}{\mathrm{d}t}\|\boldsymbol{e}_h\|_{\boldsymbol{L}^2}^2 + \nu C_\sigma \|\boldsymbol{e}_h\|_e^2 + |\boldsymbol{e}_h|_{\boldsymbol{u}_h,\mathrm{upw}}^2 \leqslant \|\partial_t \boldsymbol{\eta}\|_{\boldsymbol{L}^2}^2 + \|\boldsymbol{u}\|_{\boldsymbol{L}^\infty}\|\nabla_h \boldsymbol{\eta}\|_{\boldsymbol{L}^2}^2 \tag{49a}$$

$$+\left(1+Ch^{-2}\right)\|\nabla \boldsymbol{u}\|_{\boldsymbol{L}^\infty}\|\boldsymbol{\eta}\|_{\boldsymbol{L}^2}^2 + \left[1+\|\boldsymbol{u}\|_{\boldsymbol{L}^\infty} + C\|\nabla \boldsymbol{u}\|_{\boldsymbol{L}^\infty}\right]\|\boldsymbol{e}_h\|_{\boldsymbol{L}^2}^2. \tag{49b}$$

The essential regularity assumption $\boldsymbol{u} \in L^1\left(\boldsymbol{W}^{1,\infty}\right)$ ensures that

$$G_{\boldsymbol{u}}(t) = \int_0^t \left[1+\|\boldsymbol{u}(\tau)\|_{\boldsymbol{L}^\infty} + C\|\nabla \boldsymbol{u}(\tau)\|_{\boldsymbol{L}^\infty}\right]\mathrm{d}\tau < \infty. \tag{50}$$

Application of Gronwall's lemma [28, Lemma 6.9] together with $\boldsymbol{u}_h(0) = \boldsymbol{\pi}_s \boldsymbol{u}_0$ concludes the proof. ■

**Remark 5.7** The assumption (10) has also been used for the incompressible Euler equations ($\nu = 0$); cf. [35, 53] where $\boldsymbol{H}(\mathrm{div})$-FEM are considered. However, it has to be mentioned that (10) is very strict in case of $\nu = 0$ as there exists no inherent smoothing mechanism from the incompressible Euler operator in crosswind direction.

**Remark 5.8** In contrast to the Gronwall constants (26) and (27) for non-divergence-free FEM, the Gronwall constant (45) for divergence-free methods does not imply an explicit dependence on either $\nu^{-1}$ or any discretisation parameter (as for example the grad-div parameter $\delta$) which may involve classical Poisson locking. In this regard the results from Section 5 represent a step forwards. However, to the best of our knowledge, there does not exist numerical evidence for the sharpness of these improved estimates, thereby leaving room for further research.

**Corollary 5.9 (Velocity discretisation error convergence rate)** *Under the assumptions of the previous theorem, assume a smooth solution according to*

$$\boldsymbol{u} \in L^\infty(0,T;\boldsymbol{H}^r(\Omega)), \quad \partial_t \boldsymbol{u} \in L^2(0,T;\boldsymbol{H}^r(\Omega)), \quad r > \frac{3}{2}. \tag{51}$$

*Then, with $r_{\boldsymbol{u}} = \min\{r, k+1\}$ and a constant $C$ independent of $h$ and $\nu^{-1}$, we obtain the following convergence rate:*

$$\frac{1}{2}\|\boldsymbol{e}_h\|_{L^\infty(\boldsymbol{L}^2)}^2 + \int_0^T \left[\nu C_\sigma \|\boldsymbol{e}_h\|_e^2 + |\boldsymbol{e}_h|_{\boldsymbol{u}_h,\mathrm{upw}}^2\right]\mathrm{d}\tau \tag{52a}$$

$$\leqslant Ch^{2(r_{\boldsymbol{u}}-1)} e^{G_{\boldsymbol{u}}(T)} \int_0^T \left[h^2 |\partial_t \boldsymbol{u}|_{\boldsymbol{H}^{r_{\boldsymbol{u}}}}^2 + \left[\|\boldsymbol{u}\|_{\boldsymbol{L}^\infty} + \left(h^2 + C\right)\|\nabla \boldsymbol{u}\|_{\boldsymbol{L}^\infty}\right]|\boldsymbol{u}|_{\boldsymbol{H}^{r_{\boldsymbol{u}}}}^2\right]\mathrm{d}\tau \tag{52b}$$

**Proof:** Due to Theorem 5.3, we obtain the estimates $\|\partial_t \boldsymbol{\eta}\|_{\boldsymbol{L}^2}^2 \leqslant Ch^{2r_{\boldsymbol{u}}}|\partial_t \boldsymbol{u}|_{\boldsymbol{H}^{r_{\boldsymbol{u}}}}^2$, $\|\nabla_h \boldsymbol{\eta}\|_{\boldsymbol{L}^2}^2 \leqslant Ch^{2(r_{\boldsymbol{u}}-1)}|\boldsymbol{u}|_{\boldsymbol{H}^{r_{\boldsymbol{u}}}}^2$ and $\|\boldsymbol{\eta}\|_{\boldsymbol{L}^2}^2 \leqslant Ch^{2r_{\boldsymbol{u}}}|\boldsymbol{u}|_{\boldsymbol{H}^{r_{\boldsymbol{u}}}}^2$. The claim follows directly. ■

## 5.3 Numerical illustration of the Gronwall factor

We consider the flow of four vortices which are oppositely rotating at a fixed position in the periodic domain $\Omega = (0,1)^2$. A freely-decaying exact solution of (1) with $\boldsymbol{f} = \boldsymbol{0}$ which describes such a flow is given by

$$\boldsymbol{u}_0(\boldsymbol{x}) = \begin{bmatrix}\sin(2\pi x_1)\sin(2\pi x_2)\\ \cos(2\pi x_1)\cos(2\pi x_2)\end{bmatrix}, \quad \boldsymbol{u}(t,\boldsymbol{x}) = \boldsymbol{u}_0(\boldsymbol{x})e^{-8\pi^2\nu t}. \tag{53}$$

This example represents a generalised Beltrami flow and has already been investigated in detail, also qualitatively, in [64,65] and is called 'planar lattice flow' as well [8]. The initial velocity $\boldsymbol{u}_0$ induces a flow structure which, due to its saddle point character, is '*dynamically unstable so that small perturbations result in a very chaotic motion*' [52]. The corresponding pressure level has to be fixed; for example, by imposing the zero-mean condition. Here, we choose $\nu = 10^{-5}$ which leads to a flow where both viscous and inertia effects are present. Note that $\|\nabla \boldsymbol{u}(t)\|_{\boldsymbol{L}^\infty} = 2\pi \exp\left(-8\pi^2 \nu t\right)$.

Our aim is to demonstrate the role of the Gronwall factor for simulations over $(0,T)$ for 'large' $T$; we choose $T = 26$. This examples proves that the estimates are qualitatively sharp in the sense that the theoretically predicted exponential growth of the errors can actually be observed in practice. Note that this is not a convergence study. Related exact solutions, for example the 2D Taylor–Green problem, can also be used to show $h$ convergence at fixed (small) time instances.

All subsequent computations have been carried out using the high-order finite element library `NGSolve` [62]. The main new aspect in this work is that we use high-order FE pairs of order $k = 8$ whereas previous work in [64,65] considered only lower order methods with $k \in \{2,3\}$. Also, we now choose a different time integration procedure; namely a second-order semi-implicit BDF (SBDF2) method with constant time step size $\Delta t = 10^{-4}$; cf. [4]. The small time step makes it possible to neglect errors stemming from the time discretisation. As the implicit part we choose the Stokes-like terms (Laplacian and pressure-velocity coupling) and denote the corresponding system matrix by $M^*$. The convection part is applied explicitly and therefore the nonlinearity is shifted to the right-hand side.

We compare results on the two meshes shown in Figure 2. Note that the meshes are unstructured and therefore do not exploit the saddle-point structure of the flow. On these meshes, the $\boldsymbol{H}^1$-conforming methods under comparison are the pure Galerkin formulation of the Taylor–Hood method (Galerkin-TH8), Taylor–Hood with additional grad-div stabilisation (grad-div-TH8) with $\delta = 0.1$ (both non-divergence-free) and the divergence-free Scott–Vogelius element (SV8). The chosen $\boldsymbol{H}(\mathrm{div})$-conforming methods are based on the Brezzi–Douglas–Marini (BDM) element where one is an $\boldsymbol{H}(\mathrm{div})$-conforming DG method as in [23] (BDM8) and the other is a hybridised variant introduced in [46] (hBDM8). For the DG variant we choose $\sigma = 4k^2$ in (19) and make a corresponding choice for HDG. In terms of our analysis, both methods share the same discretisation properties but differ in computational aspects that are discussed in more detail in the Appendix.

A visualisation of the performance of the different methods can be seen in Figure 3. Let us comment on some aspects of the results. For classical Taylor–Hood elements, one observes a blow-up of the Gronwall factor due to the term $4\nu^{-1} \|\boldsymbol{u}\|^2_{L^2(\boldsymbol{L}^\infty)}$, see (26). Grad-div stabilisation with $\delta = \mathcal{O}(1)$ can considerably improve the results of the Galerkin variant. Non-div-free grad-div stabilised Taylor–Hood, div-free Scott–Vogelius FEM and div-free BDM-(H)DG show the theoretical qualitative behaviour of the exponential Gronwall factor. No immediate blow-up occurs. On the coarse mesh, $\boldsymbol{H}(\mathrm{div})$-conforming FEM provide much better results than $\boldsymbol{H}^1$-conforming FEM. In this work, we choose the (relatively) high order $k = 8$ only in order to be able to compute accurately on coarse meshes. The conclusions of the numerical experiments are in no way restricted to higher order methods. Indeed, the behaviour is consistent with the lower order case as has been observed in [65]. On the fine mesh, all *Re*-semi-robust methods perform similarly.

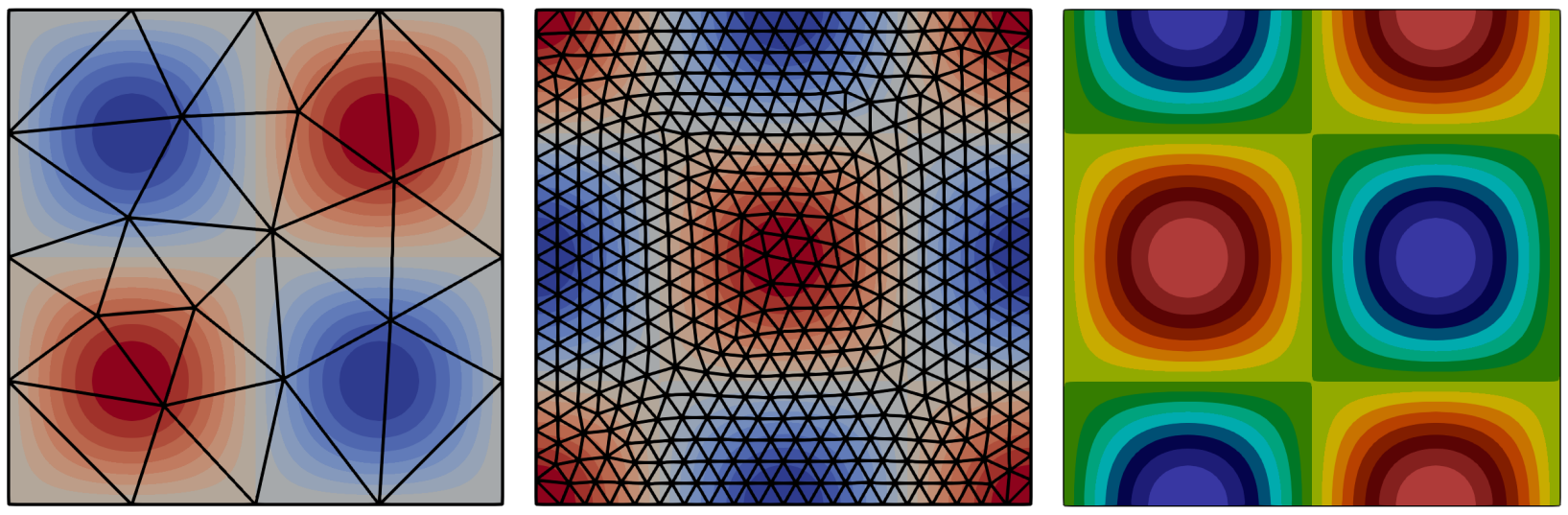

**Fig. 2** Lattice flow: Initial velocity and triangular meshes without singular vertices for the high-order FEM applied to the standing vortices problem. Left: Coarse mesh (34 triangles) with $h = 0.25$ and first component of $\boldsymbol{u}_0$; middle: fine mesh (902 triangles) with $h = 0.05$ and second component of $\boldsymbol{u}_0$; right: vorticity computed from $\boldsymbol{u}_0$.

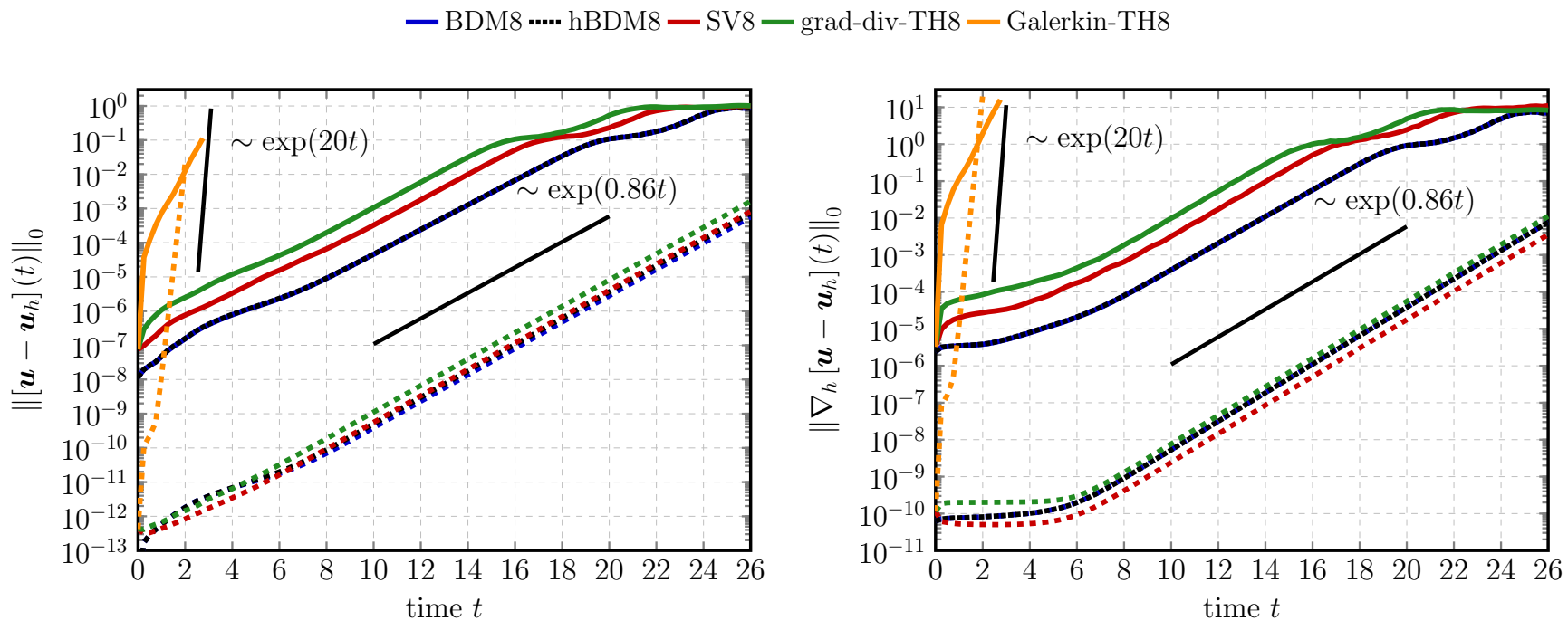

**Fig. 3** Lattice flow: Evolution of $\boldsymbol{L}^2$-norm and (broken) $\boldsymbol{H}^1$-seminorm errors for different methods. Computations on the coarse mesh are shown by solid lines whereas the fine mesh is indicated by dashed lines. The $\boldsymbol{H}(\mathrm{div})$-HDG method on the coarse mesh is shown with black dots.

## 5.4 Brief demonstration of the advantages of divergence-free methods

In this second example, we aim at showing that in certain flow configurations, divergence-free (or pressure-robust) methods outperform non-divergence-free methods immensely. To this end, consider the transient potential flow defined by $\boldsymbol{u} = \nabla\varphi$ with the harmonic potential $\varphi(t) = t\left[x^5 - 10x^3y^2 + 5xy^4\right]$. Inserting this into (1) with $\boldsymbol{f} = \boldsymbol{0}$ leads to the following exact solution; cf. [50]:

$$\boldsymbol{u} = t\begin{pmatrix} 5x^4 - 30x^2y^2 + 5y^4 \\ -20x^2y + 20xy^3 \end{pmatrix}, \quad p = -\frac{1}{2}|\nabla\varphi|^2 - \partial_t\varphi = -12.5t^2\left(x^2 + y^2\right)^4 - t^{-1}\varphi \tag{54}$$

As opposed to (1), in this example we prescribe time-dependent non-zero Dirichlet boundary conditions according to the exact solution and choose $\nu = 1$. For the solution of this problem we fix the order $k = 4$ and use the same grad-div stabilised Taylor–Hood (grad-div-TH4) and the same divergence-free $\boldsymbol{H}(\mathrm{div})$-DG (BDM4) method as for the lattice flow problem. For the time-discretisation we again employ the semi-implicit SBDF2 method with constant time-step size $\Delta t = 10^{-3}$.

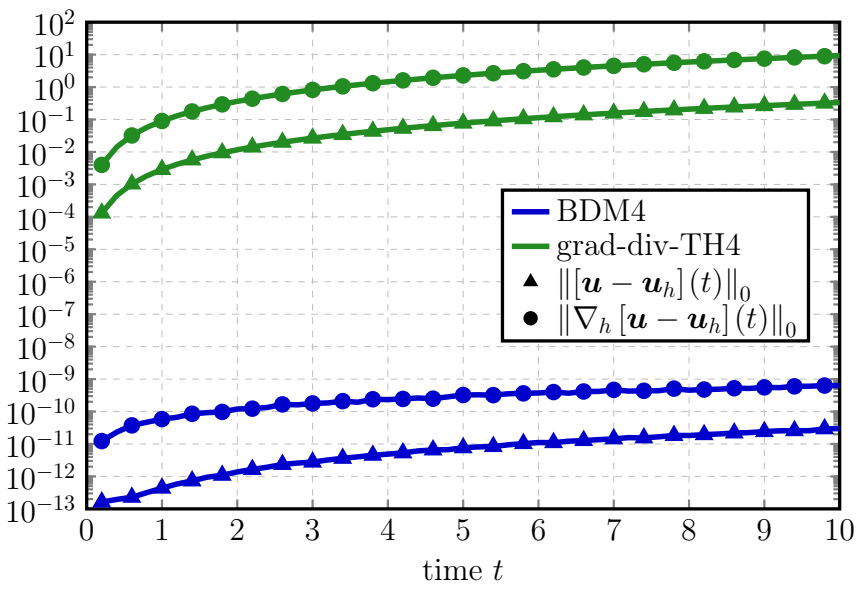

**Fig. 4** Transient potential flow: Evolution of $\boldsymbol{L}^2$-norm and (broken) $\boldsymbol{H}^1$-seminorm errors for grad-div stabilised Taylor–Hood and divergence-free $\boldsymbol{H}(\mathrm{div})$-DG. Computations are performed on an unstructured triangular mesh with $h = 0.5$ using polynomial order $k = 4$.

As can be seen in Figure 4, the divergence-free method performs about ten orders of magnitude better compared to the non-divergence-free one. In fact, the fourth-order divergence-free method delivers the exact velocity solution — up to accumulated round-off errors from time discretisation and linear solvers — in every discrete time point, while the fourth-order non-divergence-free method does not. The excellent behaviour of the pressure-robust method is somewhat surprising since the problem is discretised in space and time, simultaneously. Such a seemingly strange phenomenon can be explained by the nature of time-dependent potential flows, see [50], which are related to the elliptic Laplace problem $-\Delta\varphi(t) = 0$. Note that using $k = 4$, we have that $\boldsymbol{u} \in \boldsymbol{V}_h$ and therefore, the error which can be seen for the grad-div stabilised Taylor–Hood method originates in the lack of pressure-robustness exclusively. Indeed, the exact pressure is an eight-order polynomial and both numerical methods only use a third-order ansatz space for the discrete pressure.

## 6 Open problems

Let us comment on some open problems we deliberately circumvented in this work.

*Maximum norm estimates for $\boldsymbol{H}(\mathrm{div})$ Stokes projection:* Assumption E has been important in our analysis and yet, no rigorous mathematical proof is available. We further comment on this point in Remark 5.4.

*Sufficient conditions for the regularity assumption on $\nabla\boldsymbol{u}$:* As the regularity assumption

$$\nabla\boldsymbol{u} \in L^1(0,T;\boldsymbol{L}^\infty(\Omega)) \tag{55}$$

is crucial for a $Re$-semi-robust error analysis, sufficient conditions for this are desired. In [34, Section 7], sufficient conditions for the regularity $\nabla u \in \boldsymbol{L}^\infty$ have been derived for the stationary incompressible Navier–Stokes problem. There, in case of no-slip boundary conditions and convex polyhedral domains, a condition on the forcing term of the form $\boldsymbol{f} \in \boldsymbol{L}^{3+\varepsilon}(\Omega)$ with $0 < \varepsilon \leqslant 3/2$ is sufficient; cf. [34, Lemma 9].

*Optimality of velocity error estimate in $L^\infty(\boldsymbol{L}^2)$:* The velocity error estimate in Theorem 5.6 is optimal with respect to the spatial discretisation regarding the dissipation energy error

$$\int_0^T \nu C_\sigma \|\!|\boldsymbol{u}-\boldsymbol{u}_h|\!\|_e^2 \,\mathrm{d}\tau, \tag{56}$$

but suboptimal regarding the kinetic energy error $\frac{1}{2}\|\boldsymbol{u}-\boldsymbol{u}_h\|_{L^\infty(\boldsymbol{L}^2)}^2$. A similar result has been observed in [3,25] for inf-sup stable LPS-stabilised methods and in [33] for inf-sup stable grad-div stabilised FEM. For exactly divergence-free and $\boldsymbol{H}^1$-conforming methods, corresponding suboptimal results can be found in [64]. Let us remark that the formally quasi-optimal results for some variants of divergence-free (pressure-robust) isogeometric FEM in [29,31] are not *Re*-semi-robust. Conversely, in the equal-order case (stabilised by CIP or LPS), *Re*-semi-robust error estimates with optimal $h$-convergence rates are proved in [13, 32] but, by construction, they cannot be pressure-robust.

For exactly divergence-free $\boldsymbol{H}(\mathrm{div})$-FEM, a similar suboptimal result has been derived for the incompressible Euler problem ($\nu = 0$) in [35]. Analogously to the present work, the condition $\nabla\boldsymbol{u} \in L^1(\boldsymbol{L}^\infty)$ is crucial; see also [53]. This is a very strong regularity assumption on the solution of the Euler problem as there is no crosswind diffusion in the continuous problem. A corresponding result for the time-dependent Oseen problem ($\nu > 0$) can be found in [65], which has been extended in the present paper to problem (1) with $\nu > 0$.

For DG-FEM applied to scalar convection-diffusion problems, certain techniques can be applied to the convective term which allows an additional error order $1/2$ in case of sufficiently small viscosity $\nu \leqslant Ch$; cf. [26]. It remains an open problem whether similar techniques can be used for inf-sup stable FEM in the Navier–Stokes case. On the other hand, the suboptimality with respect to $h$ becomes less important in case of high-order FEM.

*Practically relevant boundary conditions:* In this work as well as in most of related work, the error estimates for problem (1) usually are derived under the no-slip condition $\boldsymbol{u} = \boldsymbol{0}$ for the velocity. This excludes, for example, channel-like problems with in- and outflow which are important, for example for biomedical flows (see the review [7]). Therefore, an extension of the error estimates to such more practically relevant flow problems is desired. A first attempt to grad-div stabilised FEM can be found in [2].

*Additional stabilisation and turbulence modelling:* The numerical results for the problem in Subsection 5.3 suggest that an additional stabilisation term can help to improve the results for non-divergence-free $\boldsymbol{H}^1$-FEM. Indeed, the blow-up for standard Taylor–Hood elements can be reduced dramatically by adding grad-div stabilisation. This can be explained by the fact that grad-div stabilisation can counteract problems which result from a lack of pressure-robustness; cf. [41]. Furthermore, for both Scott–Vogelius and Taylor–Hood elements, it is shown in [64,65] that an additional explicit convection stabilisation leads to some improvements as well. For $\boldsymbol{H}(\mathrm{div})$-(H)DG methods, the natural upwind mechanism takes care of dominant convection and no additional convection stabilisation is required.

However, in case of turbulent flows, an additional turbulence model, for example via subgrid-viscosity terms, might be needed. Potential candidates are a local projection stabilisation based on the Smagorinsky model [39], or residual-based eddy-viscosity methods [54]. We again want to emphasise that for $\boldsymbol{H}(\mathrm{div})$-conforming (H)DG methods with an upwind discretization for the convection, there is no need for an additional stabilisation. Therefore, the effect of explicit *turbulence modelling*, which usually has a dissipative character, can be distinguished neatly from *convection stabilisation*. This is, in our opinion, a good starting point to assess explicit turbulence models.

*Refined Gronwall estimates:* Let us return to the velocity error estimates where an exponential Gronwall factor with argument depending on $\nabla \boldsymbol{u} \in L^1(\boldsymbol{L}^\infty)$ occurs. In boundary layers, the latter term may typically depend on $\nu^{-1/2}$. It remains an open problem whether it is possible to refine the analysis, for example, based on a variational multiscale decomposition of the solution. Such an approach has been considered in [12] for two-dimensional problems (1) in case of high Reynolds numbers.

## 7 Summary and outlook

The regularity assumption $\nabla \boldsymbol{u} \in L^1(\boldsymbol{L}^\infty)$, which represents a class of flows frequently discussed in both physical and mathematical literature, leads to computable flows for exactly divergence-free FEM.

Classical inf-sup stable mixed FEM are in general not suitable for long-time integration of the time-dependent incompressible Navier–Stokes problem. This can be caused by the linear phenomenon of a lack of pressure-robustness. The nonlinear effects may reduce the computable time intervals to ultra short times.

Exactly divergence-free inf-sup stable FEM may serve as best practice examples for the time-dependent incompressible Navier–Stokes problem. In particular, the excessive growth of the exponential Gronwall factor with respect to $\nu^{-1}$ is circumvented.

Drawbacks of exactly divergence-free, inf-sup stable $\boldsymbol{H}^1$-conforming FEM stem from technical problems. Scott-Vogelius requires barycentre-refined meshes if the element order is not high-enough. Isogeometric based FEM are probably not available in standard FEM packages.

Exactly divergence-free, inf-sup stable $\boldsymbol{H}(\mathrm{div})$-conforming FEM can be constructed using Raviart–Thomas or BDM elements on arbitrary meshes. An upwind stabilisation of the convective term based on DG-FEM can be incorporated in a very natural way.

Due to their discontinuous nature, exactly divergence-free, inf-sup stable $\boldsymbol{H}(\mathrm{div})$-conforming DG-FEM can be hybridised. Such HDG-based $\boldsymbol{H}(\mathrm{div})$ methods allow highly efficient discrete solvers and, in particular, massively parallel implementation with very favourable scalability can be achieved.

Summarising, we believe that exactly divergence-free, inf-sup stable $\boldsymbol{H}(\mathrm{div})$-conforming FEM provide the most promising approach from both theoretical and practical point of view.

## Appendix: Computational aspects of $\boldsymbol{H}$(div)-conforming methods for Navier–Stokes

*Improving the efficiency of $\boldsymbol{H}(\mathrm{div})$-conforming methods:* In this section we explain how $\boldsymbol{H}(\mathrm{div})$-conforming FE methods, that are often seen as too complicated and inefficient for real application, can be made efficient. We restrict the discussion here to BDM elements as they are computationally more efficient in the context of incompressible flows compared to RT elements since they have less degrees of freedom (DOFs) for the same velocity approximation.

Choosing the pressure space $Q_h$ as the space of (discontinuous) piecewise polynomials of one degree less than the $\boldsymbol{H}(\mathrm{div})$-conforming velocity space $\boldsymbol{V}_h$ renders (16) an equality, that is, $\nabla\cdot\boldsymbol{V}_h = Q_h$. A special property of this velocity-pressure pair is that the inf-sup constant is robust in the polynomial degree leading to $hp$-optimal convergence; cf. [44] for a rigorous analysis in 2D. The strong relation $\nabla\cdot\boldsymbol{V}_h = Q_h$ can further be exploited with a smart choice of the basis functions for $\boldsymbol{V}_h$ and $Q_h$; cf. [63,67]. The *a priori* knowledge that the discrete solution will be pointwise divergence-free then allows to remove some DOFs for the velocity and all pressure unknowns except for the piecewise constants; cf. [46, Remark 1] and [45, Section 2.2.4.2]. We make use of this in our numerical experiments.

To account for the tangential discontinuity in the $\boldsymbol{H}(\mathrm{div})$-conforming FE space, a DG formulation has to be applied. This aspect can be regarded ambivalently. On the one hand, the discontinuous nature of the tangential component offers the possibility of applying an upwind discretisation for the convection, cf. (22), which results in stable discretisations also in the convective limit [35] without adding too much dissipation compared to most convection stabilisations of $\boldsymbol{H}^1$-FEM. On the other hand, the DG formulation results in computationally less attractive features. Due to the break-up of the tangential continuity, several DOFs for the velocity are multiplied compared to $\boldsymbol{H}^1$-conforming methods. Even worse, the number of couplings in a corresponding system matrix increases which results in much higher computational costs for (direct and iterative) solvers of linear systems.

Several measures can be taken to compensate for these costs. To this end, we briefly discuss the concept of hybridisation in the context of DG methods [22]. To reduce the couplings of neighbouring elements, additional unknowns on the facets are introduced (which typically approximate the trace of the unknown field). These additional unknowns are used to replace the direct couplings of neighbouring elements with couplings between element unknowns and the facet unknowns. Due to the lower dimension of the facets, this reduces the overall amount of couplings especially in the higher order case. More importantly, it allows for static condensation, i.e. the elimination of interior unknowns by a local Schur complement strategy which reduces the number of DOFs for which a global linear system needs to be solved.

Depending on the problem at hand there are many ways to make use of hybridisation. For an overview we refer to the review article [18]. For Stokes and Navier–Stokes discretisations many variants have been considered; see, for instance, [21,19,20]. Exactly divergence-free HDG methods have also been considered in [17] and [57,58] where additional facet unknowns can be used to enforce normal continuity on a standard DG space which circumvents the construction of $\boldsymbol{H}(\mathrm{div})$-conforming FE spaces. Here, we use the formulation presented in [46] where, additionally to an $\boldsymbol{H}(\mathrm{div})$-conforming FE space $\boldsymbol{V}_h$ for the velocity and a discontinuous pressure space $Q_h$, facet unknowns are introduced only for the tangential component of the velocity. The DG terms in the variational formulation are then adjusted correspondingly. Finally, the element unknowns of the $\boldsymbol{H}(\mathrm{div})$-conforming FE space couple with neighbour elements only through facet unknowns. These facet unknowns are either the DOFs for the normal continuity of $\boldsymbol{V}_h$ or the additional facet unknowns. All remaining velocity unknowns, as well as the pressure unknowns, have only element local couplings such that these—except for the mean element pressure—can be eliminated during static condensation; cf. Figure 5 for a sketch.

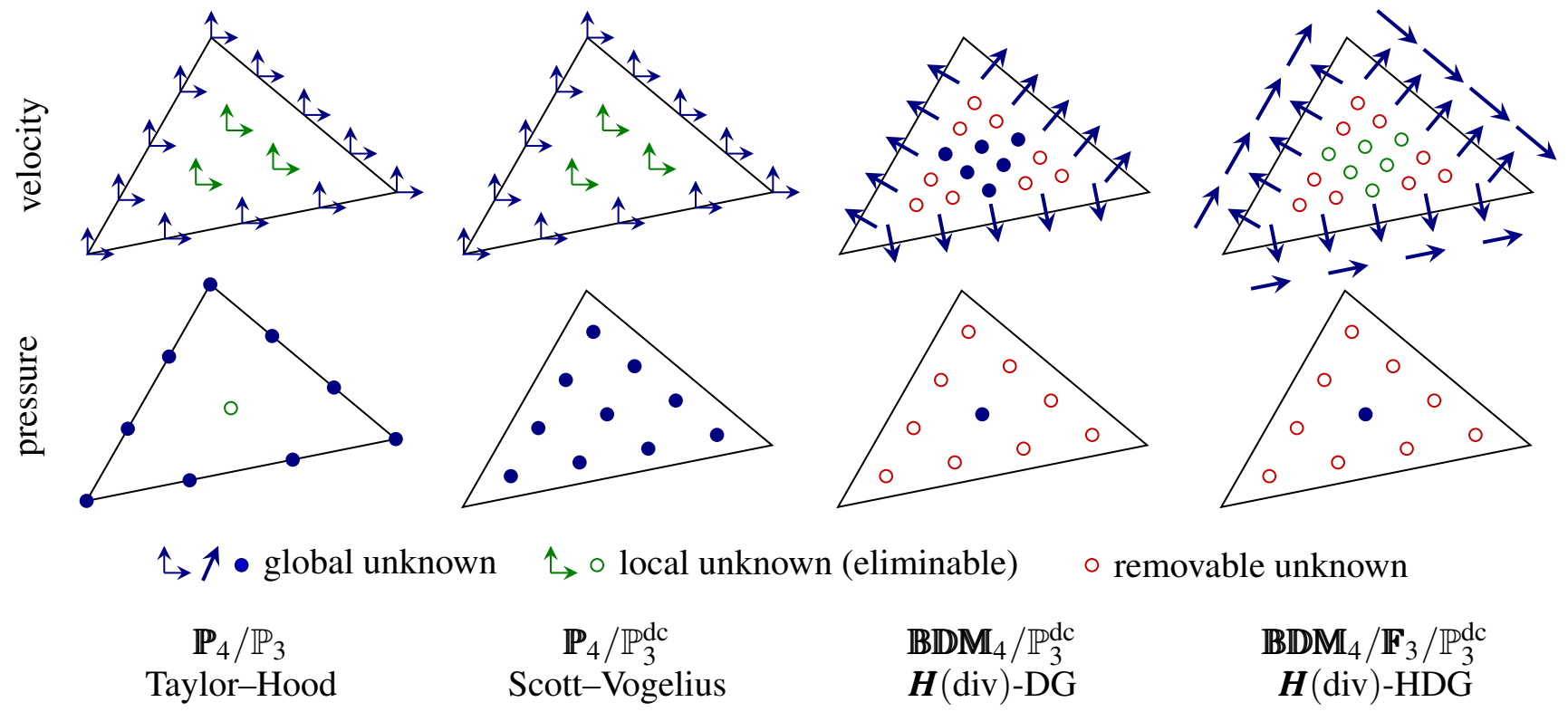

**Fig. 5** Sketch of fourth order FE discretisations with different types of unknowns for velocity and pressure: unknowns that can be remove beforehand if a suitable basis is used, local unknowns that can be eliminated by static condensation and the remaining global unknowns. The space $\mathbb{F}_3$ in the $\boldsymbol{H}(\text{div})$-HDG method is the space of vector-valued functions that are tangential polynomials up to degree three on every facet.

In the viscosity dominated case hybridisation can be optimised further so that only facet unknowns of one degree less need to be considered; cf. [46,55,56]. A similar optimisation can also be made for the unknowns for the normal continuity by relaxing the $\boldsymbol{H}(\text{div})$-conformity slightly. We do not treat this here but instead refer to [42]. To make use of these *superconvergence* properties of HDG methods we apply—as suggested in [46]—an operator splitting time integration method where the convection operator is treated only explicitly while the remaining time-independent operators are treated implicitly. Note that such an operator splitting is not only desirable for hybrid DG discretisations. Several time integration methods allow for such a splitting; cf. [46, Section 3]. For the experiments in Section 5.3 a second-order implicit-explicit BDF2 method has been used.

*Some performance comparisons for the numerical study in Section 5.3:* In Section 5.3 the errors for Taylor–Hood, Scott–Vogelius, BDM and the hybridised BDM FE discretisation on two different meshes are compared. At this point, this study shall be complemented with information on the computational costs of the methods. The results are shown in Table 1 where we make this comparison only in terms of the following measures. Firstly, the numbers of DOFs for velocity and pressure ($\#\{\boldsymbol{u}\,\text{DOFs}\}$, $\#\{p\,\text{DOFs}\}$, $\#\{\text{DOFs}\}$) are compared. Secondly, we consider the same numbers that remain in a global linear system after static condensation and a potential reduction of the basis (in brackets). Thirdly, the non-zero entries in the global matrix $M^*$ before ($\#\{\text{nz}(M^*)\}$) and after reduction and static condensation (in brackets) are considered. Note that these numbers can only give an indication of the computational efficiency of the methods. Many different practically relevant aspects, as for example parallelisability or the availability and performance of suitable preconditioners, are not reflected in these numbers.

**Table 1** Overview of meshes, DOFs and non-zero entries of $M^*$. Abbreviations of different methods with $k=8$: Non-div-free $\boldsymbol{H}^1$ Taylor–Hood (TH8), div-free $\boldsymbol{H}^1$ Scott–Vogelius (SV8) and div-free $\boldsymbol{H}(\mathrm{div})$ Brezzi–Douglas–Marini (BDM8) together with the hybridised variant (hBDM8). In brackets are the numbers after reduction of the basis and static condensation.

| Mesh | Method | $\#\{\boldsymbol{u}\,\text{DOFs}\}$ | $\#\{p\,\text{DOFs}\}$ | $\#\{\text{DOFs}\}$ | $\#\{\text{nz}(M^*)\}$ |
|---|---|---|---|---|---|
| Coarse ($h=0.25$) 34 triangles | TH8 | 2306 (748) | 890 (323) | 3196 (1071) | 465K (128K) |
| | SV8 | 2306 (748) | 1224 (1224) | 3530 (1972) | 479K (223K) |
| | BDM8 | 2673 (1483) | 1224 (34) | 3897 (1517) | 1.93M (327K) |
| | hBDM8 | 3204 (867) | 1224 (34) | 4428 (901) | 779K (77.2K) |
| Fine ($h=0.05$) 902 triangles | TH8 | 58 370 (19 844) | 22 380 (8569) | 80 750 (28 413) | 12.3M (3.38M) |
| | SV8 | 58 370 (19 844) | 32 472 (32 472) | 90 842 (52 316) | 12.7M (5.93M) |
| | BDM8 | 69 363 (37 793) | 32 472 (902) | 101 835 (38 695) | 51.3M (8.70M) |
| | hBDM8 | 81 900 (23 001) | 32 472 (902) | 38 695 (23 903) | 20.7M (2.05M) |

Regarding static condensation in the Taylor–Hood method, independent of the grad-div stabilisation, we can eliminate all interior unknowns for velocity and pressure. On general meshes, the pressure unknowns for the Scott–Vogelius element cannot be eliminated and hence, static condensation is only applied with respect to the interior velocity DOFs. We note that on barycentre refined meshes static condensation can also be applied for the pressure unknowns; cf. [24]. In case of a DG formulation with BDM elements we utilise the special basis introduced in [63,67] to eliminate some velocity unknowns and all pressure unknowns except for the mean element pressure. However, static condensation cannot be applied to any additional DOFs due to the DG couplings. Note that this could potentially be improved slightly by choosing a nodal basis similar to the one in [37] where interior unknowns only couple with the boundary nodes of neighbouring elements. To the best of the authors' knowledge, such a basis has not yet been proposed for an $\boldsymbol{H}(\mathrm{div})$-conforming FE space. For the hybridised DG method we can apply the reduction of the basis for the $\boldsymbol{H}(\mathrm{div})$-conforming FE space as well as static condensation. Note that in this work, the formulation from [46] is used which only involves tangential facet unknowns of degree 7. The results are shown in Table 1.

We observe that the effect of the basis reduction and especially the hybridisation reduces the computational costs of the $\boldsymbol{H}(\mathrm{div})$-conforming methods drastically, thereby rendering them competitive not only in terms of accuracy, cf. Section 5.3, but also in terms of computing time; see also the benchmark results in [46, Section 4.5].